\documentclass[12pt]{article} 
\usepackage{amssymb,amsmath,amsthm,tikz,multirow,graphicx,xcolor,subfigure}
\usetikzlibrary{calc,arrows}

\title{Tilings of Flat Tori by Congruent Hexagons}
\author{Xinlu Yu, Erxiao Wang\thanks{Corresponding author (wang.eric@zjnu.edu.cn).  Research was supported by National Natural Science Foundation of China NSFC-RGC 12361161603 and Key Projects of Zhejiang Natural Science Foundation LZ22A010003.}, 
	Zhejiang Normal University \\
	Min Yan\thanks{Research was supported by NSFC-RGC Joint Research Scheme N-HKUST607/23 and Hong Kong RGC General Research Fund 16305920.}, 
	Hong Kong University of Science and Technology}

\usepackage[hidelinks, hyperindex]{hyperref}

\hyperref[sec:function]{}

\definecolor{emerald}{rgb}{0.31, 0.78, 0.47}

\newcommand{\sub}{\subset}

\newcommand{\mc}{\mathcal}
\newcommand{\bb}{\mathbb}

\newtheorem{theorem}{Theorem}
\newtheorem{lemma}[theorem]{Lemma}

\newtheorem*{theorem*}{Theorem}

\theoremstyle{definition}
\newtheorem*{definition*}{Definition}
\newtheorem*{example*}{Example}
\newtheorem*{case*}{Case}
\newtheorem*{subcase*}{Subcase}

\theoremstyle{remark}

\numberwithin{equation}{section}

\begin{document}
	
\date{}
\maketitle

\begin{abstract}
Convex hexagons that can tile the plane have been classified into three types. For the generic cases (not necessarily convex) of the three types and two other special cases, we classify tilings of the plane under the assumption that all vertices have degree $3$. Then we use the classification to describe the corresponding hexagonal tilings of flat tori and their moduli spaces.
	
\end{abstract}

Reinhardt \cite{re} classified convex hexagons that can tile the Euclidean plane into three types (also see \cite{bo,zo1}). However, it is not known whether the classification is still valid for non-convex hexagons. Moreover, the way the hexagon tiles the plane is not well understood. In fact, it may still be an open problem whether there are anisohedral hexagons   \cite{gs,bmp}, although the method of Rao \cite{ra} can conceivably be used to solve the problem. 

It is well known that, under some mild condition on the size of hexagons, a (not necessarily monohedral) hexagonal planar tiling has arbitrarily large patches in which all vertices have degree three \cite{re,bo}. Therefore we consider monohedral hexagonal planar tilings, such that all vertices have degree three. For three generic cases of Reinhardt's hexagons, and the generic centrally symmetric hexagon, we show that such planar tiling is unique and periodic (Theorems \ref{combo1} to \ref{combo4}). For another case, the tiling can be more flexible, and we give complete description (Theorem \ref{combo5}). 

The universal cover of a torus tiling is a periodic planar tiling. Consequently, we know tilings of the torus by the hexagons in Theorems \ref{combo1} to \ref{combo5}. Since a finite covering of a torus tiling is still a torus tiling, the description of torus tilings is reduced to the coverings of ``minimal'' torus tilings (which are not proper coverings of other torus tilings). In Theorem \ref{cover}, we describe all coverings in terms of the Hermite normal form. 

For the hexagons in Theorems \ref{combo1} to \ref{combo4}, we identify the unique minimal tilings and find that they  allow some free parameters. Then we describe the moduli spaces of all the minimal tilings, similar to our earlier study of the moduli spaces of pentagonal subdivision tilings \cite{lwy}.

\section{Basic Facts and Techniques}

For the basics of tilings, we adopt the definitions and conventions in the books by Gr\"{u}nbaum and Shephard \cite{gs} and by Adams \cite{ad}. The vertices and edges of a polygon are renamed \emph{corners} and \emph{sides}. A \emph{tiling} of a plane or surface is a partition into closed polygons, called \emph{tiles}. A {\em vertex} is the meeting place of at least three tiles. An {\em edge} is an arc in the intersection of two tiles, such that the two ends are vertices, and there are no vertices in the interior of the arc. 

\begin{figure}[htp]
\centering
\begin{tikzpicture}[>=latex,scale=1]
	
\draw
	(-0.8,-0.8) -- (-1.6,-0.8) -- (-1.6,0) -- (1.6,0) -- (2.4,0.8)
	(0,0) -- (-0.8,-0.8) -- (0,-1.6) -- (0.8,-0.8) -- (0,0)
	(-0.8,0) -- (-0.8,0.8) -- (0.8,0.8) -- (2.4,-0.8) -- (0.8,-0.8)
	(1.6,0) -- (0.8,0.8)
	(0.8,0.8) -- (1.6,1.6) -- (3.2,0) -- (2.4,-0.8);

\fill
	(1.6,0) circle (0.05);
\filldraw[fill=white]
	(0,0) circle (0.05);

\node[inner sep=0.5, draw, shape=circle] at (0,0.4) {\tiny 1};
		
\end{tikzpicture}
\caption{Full vertex $\bullet$ and half vertex $\circ$.}
\label{vertex}
\end{figure}
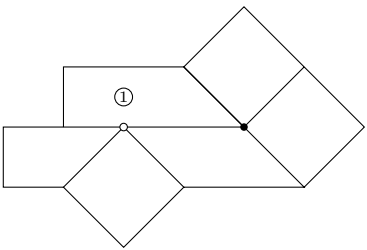

We emphasize that corners and sides refer to polygons or tiles, and vertices and edges refer to the tiling. If the vertices and edges coincide with the corners and sides, then we say the tiling is \emph{side-to-side}\footnote{This is usually called {\em edge-to-edge} in the literature.}. If a tiling is not side-to-side, then a vertex may lie in the interior of a side, such as the vertex $\circ$ in Figure \ref{vertex}. Then we call the vertex a {\em half vertex}. Otherwise the vertex is a {\em full vertex}, such as the vertex $\bullet$ in Figure \ref{vertex}. A vertex is a full vertex if and only if it is the corner of all the tiles at the vertex. 

The {\em degree} of a vertex is the number of corners at the vertex. For example, the full vertex $\bullet$ in Figure \ref{vertex} has degree $4$, and the half vertex $\circ$ has degree $3$. In this paper, we only consider tilings in which all sides are straight, and full vertices have degrees $\ge 3$, and half vertices have degree $\ge 2$. Such a degree condition is equivalent to the definition of {\em proper normal polygonal tiling} in \cite{gs}.

The assumption of straight sides implies that any half vertex is inside the interior of only one side, and the property is used only for the proof of Lemma \ref{deg3}. Therefore the straight side assumption can be replaced by the property of any half vertex being inside the interior of only one side. The sides need not be straight.

\subsection{Vertex}

The following is the simpler version of Theorem 1 of \cite{bo}. We also give a simpler proof.

\begin{lemma}\label{deg3}
A hexagonal tiling of the torus (or Klein bottle) is side-to-side, and all vertices are full vertices of degree $3$. 
\end{lemma}

\begin{proof}
In a hexagonal tiling of the torus, let $v,h,e,f$ be the numbers of full vertices, half vertices, edges, and tiles. Then we have
\begin{align}
(v+h)-e+f &=0, \label{eq1} \\
6f+h &=2e. \label{eq2}
\end{align}
The first equality is the Euler number of the torus (or Klein bottle). For the second equality, we note that the number of edges in a tile is $6$ plus the number of half vertices on the sides of the tile. Since each half vertex has only one side, corresponding to only one tile, the sum $\sum$ of all these numbers is $6f+h$. On the other hand, since each edge is shared by exactly two tiles, the sum $\sum$ is $2e$. 

Canceling $e$ from the equations \eqref{eq1} and \eqref{eq2}, we get $2v+h=4f$. 

Let $v_k$ be the number of full vertices of degree $k$. Let $h_l$ be the number of half vertices of degree $l$. Then we have
\begin{align*}
v &=\sum_{k\ge 3}v_k=v_3+v_4+v_5+\cdots, \\
h &=\sum_{l\ge 2}h_l=h_2+h_3+h_4+\cdots, \\
6f &=\sum_{k\ge 3}kv_k+\sum_{l\ge 2}lh_l.
\end{align*}
The first two equalities follow directly from the definition. The last equality is due to two ways of counting the total number of corners in a tiling. Substituting the three equalities into $2v+h=4f$, we get 
\[
2v_4+4v_5+6v_6+\cdots+h_2+3h_3+5h_4+\cdots=0.
\]
This means no full vertices of degree $\ge 4$, and no half vertices. 
\end{proof}

\subsection{Flat Tori}

For the congruence of polygons to make sense in a torus, the torus must have a flat metric. The universal cover of a flat torus is the Euclidean plane ${\bb C}$ with the usual metric, and the flat torus is the quotient of the Euclidean plane by a lattice $\Lambda={\bb Z}\alpha+{\bb Z}\beta\sub {\bb C}$, where $\alpha,\beta$ form a real basis of ${\bb C}$. The congruence in the torus comes from the usual congruence in the plane ${\bb C}$. 

A {\em scaling} of the torus corresponds to multiplying a positive number $c$ to the lattice. We consider a metric and its scaling (also called homothety) to be equivalent. In other words, we regard the tori corresponding to $\Lambda$ and $c\Lambda={\bb Z}c\alpha+{\bb Z}c\beta$ ($c>0$) to be isometric. 

Two tilings ${\mc T}_1$ and ${\mc T}_2$ on (flat) tori $T_1$ and $T_2$ are {\em equivalent} if there is an isometry (up to scaling) $T_1\to T_2$ that takes ${\mc T}_1$ to ${\mc T}_2$. Moreover, the covering of a torus tiling is still a torus tiling (the two tori may not be isometric). Therefore we need to consider isometric classes of tori, the equivalence of tilings on a fixed torus by self isometries, and the coverings of torus tilings.

The isometric classes of tori is part of the Teichm\"uller theory. We give a brief outline relevant to the torus tiling. By exchanging $\alpha,\beta$ (which does not change $\Lambda$) if necessary, we may assume $\text{Im}\frac{\beta}{\alpha}>0$. Geometrically, this means $\alpha\to\beta$ is in counterclockwise direction. Then we may multiply the lattice by $\alpha^{-1}$ (which is rotation and scaling, i.e., an isometry) to get the standard presentation $\Lambda_{\tau}={\bb Z}1+{\bb Z}\tau$ of the lattice. Here $\tau=\frac{\beta}{\alpha}$ lies in the open upper half plane $\mathbb{H}^{2}=\{\tau\in {\bb C}\colon \text{Im}\tau>0\}$. 

An isometry between two tori $T_{\tau}={\bb C}/\Lambda_{\tau}$ and $T_{\tau'}={\bb C}/\Lambda_{\tau'}$ means $\Lambda_{\tau}$ and $\Lambda_{\tau'}$ are isometric: $\Lambda_{\tau'}=\lambda\Lambda_{\tau}$ for some $\lambda\in {\bb C}-0$. In other words, we have the relation 
\[
\begin{pmatrix}
1 \\ \tau'
\end{pmatrix}
=
\lambda\begin{pmatrix}
a & b \\ c & d
\end{pmatrix}
\begin{pmatrix}
1 \\ \tau
\end{pmatrix},\quad
\begin{pmatrix}
a & b \\ c & d
\end{pmatrix}\in SL_2(\mathbb{Z}).
\]
Therefore $\tau$ and $\tau'$ are related by the action of $SL_2({\bb Z})$ on ${\mathbb H}^2$
\[
\tau'=\mu(\tau)=\frac{c+d\tau}{a+b\tau},\quad
\mu=\begin{pmatrix}
a & b \\ c & d
\end{pmatrix}\in SL_2({\bb Z}).
\]
The moduli space of isometric classes of tori is $\mathbb{H}^{2}/SL_2(\mathbb{Z})$.

\subsection{Coverings of Tori}

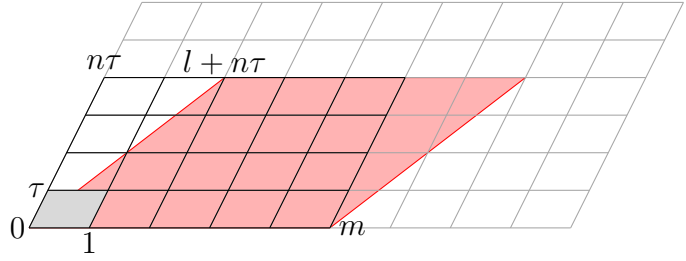
\begin{figure}[htp]
\centering
\begin{tikzpicture}[>=latex,scale=1]

\filldraw[draw=red, fill=red!30]
	(0,0) -- (4,0) -- (6.6,2) -- (2.6,2) -- (0,0); 

\fill[gray!30]
	(0,0) -- (0.8,0) -- (1.05,0.5) -- (0.25,0.5);

\foreach \a in {0,...,9}
\draw[gray!70]
	(0.8*\a,0) -- ++(1.5,3);

\foreach \a in {0,...,6}
\draw[gray!70]
	(0.25*\a,0.5*\a) -- ++(7.2,0);
				
\foreach \a in {0,...,5}
\draw
	(0.8*\a,0) -- ++(1,2);

\foreach \a in {0,...,4}
\draw
	(0.25*\a,0.5*\a) -- ++(4,0);	

\node at (0.1,0.5) {$\tau$};
\node at (-0.15,0) {$0$};
\node at (0.8,-0.2) {$1$};
\node at (4.3,0) {$m$};
\node at (1,2.2) {$n\tau$};
\node at (2.6,2.2) {$l+n\tau$};

\end{tikzpicture}
\caption{Covering of ${\bb C}/\Lambda_{\tau}$ labeled by $(m,n;l)$, with $m,n\in {\bb Z}_+$ and $0\le l<m$.}
\label{covering}
\end{figure}

A covering of the torus $T_{\tau}={\bb C}/\Lambda_{\tau}$ is a torus ${\bb C}/\Lambda$ with the lattice satisfying $\Lambda\sub \Lambda_{\tau}$ (called \emph{sub-lattice}). In Figure \ref{covering}, the lattice $\Lambda_{\tau}$ is described by the black and gray lines, with the fundamental domain being the parallelogram spanned by $1$ and $\tau$. 

The lattice $\Lambda$ is generated by a real basis $a+b\tau,c+d\tau\in \Lambda_{\tau}$, with $a,b,c,d\in {\bb Z}$ satisfying $ad-bc\ne 0$. The basis can be changed while still fixing $\Lambda$, such as the exchange of the two vectors, or
\[
\Lambda
={\bb Z}(a+b\tau)+{\bb Z}(c+d\tau)
={\bb Z}(a+b\tau)+{\bb Z}((c-ka)+(d-kb)\tau)
=\cdots.
\]
Therefore we may apply the Euclidean algorithm to the coefficients $a,b,c,d$, and improve the basis to $m,l+n\tau\in \Lambda$, where $n$ is the greatest common divisor of $b$ and $d$, and $m,n$ are positive integers. We further note that $l+n\tau$ can be modified by adding or subtracting multiples of $m$. Therefore we may assume $0\le l<m$. Then we get the unique expression, called the {\em Hermite normal form} \cite{he}
\[
\Lambda=\Lambda_{\tau}(m,n;l)={\bb Z}m+{\bb Z}(l+n\tau),\quad
m,n\in {\bb Z}_+,\; 
0\le l<m.
\]
The covering of $T_{\tau}$ is ${\bb C}/\Lambda_{\tau}(m,n;l)=T_{\frac{l+n\tau}{m}}$.

A tiling of $T_{\tau}$ is given by some tiles in the {\em fundamental domain} of $\Lambda_{\tau}$ (the gray parallelogram in Figure \ref{covering}). Then the tiles are repeated throughout the lattice $\Lambda_{\tau}$. Those tiles in the fundamental domain of $\Lambda_{\tau}(m,n;l)$ (the red parallelogram) form the covering tiling of the covering torus $T_{\frac{l+n\tau}{m}}$. In particular, for any rational numbers $r,s>0$, the same prototile tiles $T_{r+s\tau}$. Therefore the tori that can be tiled by the prototile form a dense subset of the moduli space $\mathbb{H}^{2}/SL_2(\mathbb{Z})$ of flat tori. 

The covering of torus tiling induces the concept of minimal tiling. A torus tiling ${\mc T}$ is {\em minimal}, if ${\mc T}$ covering another torus tiling ${\mc T}'$ implies ${\mc T}={\mc T}'$. Since any torus tiling consists of only finitely many tiles, it covers a minimal torus tiling in the manner described above. Therefore the classification of all torus tilings is reduced to the classification of minimal torus tilings. 

\begin{theorem}\label{cover}
Given a prototile, a torus tiling by the prototile covers a minimal one. Moreover, all the tori that can be tiled by the prototile form a dense subset of the moduli space of flat tori.
\end{theorem}

We remark that the discussion about the covering tiling does not actually use hexagons. In fact, the theorem remains valid for any finite protoset.

\section{Tilings by Generic Reinhardt Hexagons}

\begin{figure}[htp]
	\centering
	\begin{tikzpicture}[>=latex,scale=1]
		
		% label
		
		\foreach \a in {0,...,5}
		{
			\draw[rotate=60*\a]
			(0:1.2) -- (60:1.2);
			
			\node at (-60+60*\a:1) {\footnotesize \a};
			\node at (60*\a-30:0.85) {\footnotesize $\bar{\a}$};
		}
		
		% type 1
		
		\begin{scope}[xshift=3cm]
			
			\draw
			(-60:1.2) -- (0:1.2) -- (60:1.2)
			(120:1.2) -- (180:1.2) -- (240:1.2);
			
			\draw[red]
			(120:1.2) -- (60:1.2)
			(-120:1.2) -- (-60:1.2);
			
			\draw[->]
			(90:0.3) -- (90:1);
			\draw[->]
			(-90:0.3) -- (-90:1);
			
			\node at (0,0) {\small parallel};
			
			\node at (0,-1.3) {\small Type I};
			
		\end{scope}
		
		% type 2
		
		\begin{scope}[xshift=6cm]
			
			\draw
			(0:1.2) -- (-60:1.2)
			(180:1.2) -- (240:1.2);
			
			\draw[red]
			(120:1.2) -- (180:1.2)
			(60:1.2) -- (0:1.2);
			
			\draw[blue]
			(120:1.2) -- (60:1.2)
			(-120:1.2) -- (-60:1.2);
			
			\draw[->]
			(0:0.4) -- (0:1.1);
			\draw[->]
			(120:0.3) -- (120:1.1);
			\draw[->]
			(-60:0.3) -- (-60:1.1);
			
			\node at (-0.3,0) {\small $\Sigma=2\pi$};
			
			\node at (0,-1.3) {\small Type II};
			
		\end{scope}
		
		% type 3
		
		\begin{scope}[xshift=9cm]
			
\draw[emerald]
	(-60:1.2) -- (-120:1.2) -- (180:1.2);

\draw[red]
	(-60:1.2) -- (0:1.2) -- (60:1.2);

\draw[blue]
	(60:1.2) -- (120:1.2) -- (180:1.2);

\draw[->]
	(0:0.3) -- (0:1.1);
\draw[->]
	(120:0.35) -- (120:1.1);
\draw[->]
	(240:0.35) -- (240:1.1);
			
			\node at (0,0) {\small $\frac{2}{3}\pi$};
			
			\node at (0,-1.3) {\small Type III};
			
		\end{scope}
		
	\end{tikzpicture}
	\caption{Labels for corners and sides, and three types of prototiles.}
	\label{notation}
\end{figure}
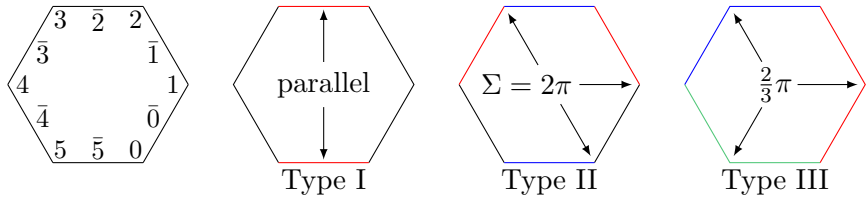

The first of Figure \ref{notation} shows a hexagonal prototile. We label the corners by $i\in {\bb Z}_6$, and denote the angle value of the corner $i$ by $[i]$. We also label the side connecting $i$ and $i+1$ by $\bar{i}$, and denote the length of the side by $|\bar{i}|$. By Lemma \ref{deg3}, in a hexagonal tiling of a torus, all sides are edges.

By Lemma \ref{deg3}, the universal cover of a tiling on a flat torus is a side-to-side tiling of the plane, such that all vertices have degree $3$. Therefore a hexagonal prototile that tiles a torus also tiles the plane. By Reinhardt \cite{re}, if the hexagon is convex, then it is one of the three types:
\begin{itemize}
	\item Type I: $[0]+[1]+[2]=[3]+[4]+[5]=2\pi$, $|\bar{2}|=|\bar{5}|$.
	\item Type II: $[0]+[1]+[3]=[2]+[4]+[5]=2\pi$, $|\bar{1}|=|\bar{3}|$, $|\bar{2}|=|\bar{5}|$.
	\item Type III: $[1]=[3]=[5]=\tfrac{2}{3}\pi$, $|\bar{0}|=|\bar{1}|$, $|\bar{2}|=|\bar{3}|$, $|\bar{4}|=|\bar{5}|$.
\end{itemize}
The three types are illustrated in Figure \ref{notation}. The red, blue, and green mean sides of the same lengths. The normal sides can have any lengths.

Under extra generic conditions and all vertices having degree $3$, we classify side-to-side tilings of the plane by the three types. We find they are all isohedral and periodic. Then we identify the minimal torus tilings. Recall that Gr\"{u}nbaum and Shephard classified and listed $13$ hexagonal isohedral types in page $481$ of \cite{gs}. For the purpose of this paper, we describe the symmetry groups of the three generic tilings in detail. 

In the subsequent discussion, the tiles in a tiling will be labeled by circled numbers \raisebox{.2pt}{\textcircled{\raisebox{0.3pt} {\scriptsize $j$}}} in pictures, and denoted $t_j$. We add subscripts to indicate which tile certain corners and sides belong to. Therefore $i_j$ and $\bar{i}_j$ are the corner $i$ and side $\bar{i}$ in the tile $t_j$. We introduced these notations in our earlier work on spherical tilings \cite{wy1,wy2,awy}.

\subsection{Generic Type I Tiling}

\begin{theorem}\label{combo1}
A type ${\rm I}$ hexagonal prototile satisfies 
\[
[0]+[1]+[2]=[3]+[4]+[5]=2\pi,\quad
|\bar{2}|=|\bar{5}|.
\]
Suppose it further satisfies the following generic conditions:
\begin{itemize}
\item Each of $\bar{0},\bar{1}$ has a different length from the other five sides.
\item $|\bar{3}|\ne |\bar{4}|$.
\end{itemize}
Then its tiling of the plane, such that all vertices have degree $3$, is uniquely given by Figures \ref{tiling1type} and \ref{tiling1a}. The tiling is isohedral, periodic, and covers a minimal torus tiling with two tiles. 
\end{theorem}

\begin{figure}[htp]
\centering
\begin{tikzpicture}[>=latex,scale=1]

\draw[gray!70]
	(0,0) -- (4,0) -- (5,2) -- (1,2) -- (0,0); 	

\draw
	(1,-0.2) -- (0,0) -- (-0.2,0.8) -- (1,2) -- (2,1.8) -- (2.8,1) -- (1,-0.2) -- (1.8,-1) -- (2.8,-1.2) -- (4,0) -- (3.8,0.8);

\draw[->, thick]
	(3.8,0.8) -- ++(-1,0.2);

\node at (4.2,0) {$1$};
\node at (1,2.2) {$\tau$};
\node at (3.4,1.05) {$\sigma$};

\fill
	(0,0) circle (0.05)
	(1,2) circle (0.05)
	(2.8,-1.2) circle (0.05)
	(3.8,0.8) circle (0.05);

\begin{scope}[red]

\draw
	(0.8,1) -- (0,0) -- (-0.2,2.4) -- (1,2) -- (1.8,3) -- (3,1.4) -- (0.8,1) -- (2,-0.6) -- (2.8,0.4) -- (4,0) -- (3.8,2.4);

\draw[->, thick]
	(3.8,2.4) -- (3,1.4);

\node at (3.35,2.05) {$\sigma$};

\end{scope}

\node[inner sep=0.5, draw, fill=white, shape=circle] at (1.2,1) {\tiny 1};
\node[inner sep=0.5, draw, fill=white, shape=circle] at (2.6,-0.2) {\tiny 2};
\node[inner sep=0.5, draw, fill=white, shape=circle] at (3.6,1.8) {\tiny 3};

\begin{scope}[xshift=6cm]

\fill[gray!30]
	(3.8,0.8) -- (-0.2,0.8) -- (1,2);
\fill[red!30]
	(3.8,0.8) -- (-0.2,0.8) -- (0,0);
\fill[color=blue!30]
	(3.4,3.4) -- (2.4,1.4) -- (3.8,0.8) -- (4.8,2.8) -- (5.1,3.4)
	(0.5,-2.4) -- (1.9,0.4) -- (3.8,0.8) -- (2.2,-2.4);
\fill[emerald!30]
	(2.8,-1.2) -- (4,0) -- (3.8,0.8) -- (5,2) -- (4.8,2.8);

\draw[gray]
	(0,0) -- (4,0) -- (5,2) -- (1,2) -- (0,0);

\begin{scope}[red]

\draw
	(1.4,0.5) -- (2.4,0.3) -- (0.4,-1.5);

\draw[->, thick]
	(3.8,0.8) -- (1.4,0.5);

\draw
	(2.8,-1.2) -- ++(-2.4,-0.3)
	(0,0) -- (2.4,0.3);

\end{scope}

\begin{scope}[blue]

\draw
	(3.6,2.8) -- (0.2,-2) -- (2.6,0.8)
	(2.8,-1.2) -- ++(-0.2,2)
	(0,0) -- (0.2,-2);;

\draw[->, thick]
	(3.8,0.8) -- ++(-0.2,2);

\end{scope}

\begin{scope}[emerald]

\draw
	(4.6,1.8) -- (-0.8,-1) -- (3.6,-0.2)
	(2.8,-1.2) -- ++(0.8,1)
	(0,0) -- (-0.8,-1);

\draw[->, thick]
	(3.8,0.8) -- (4.6,1.8);

\end{scope}

\draw
	(2.8,-1.2) -- (4,0) -- (3.8,0.8) -- (5,2) -- (4.8,2.8)
	(0,0) -- (-0.2,0.8) -- (1,2)
	(1.8,-1) -- (1,-0.2) -- (2.8,1) -- (2,1.8)
	(0,0) -- ++(1,-0.2)
	(1,2) -- ++(1,-0.2)
	(2.8,-1.2) -- ++(-1,0.2);

\draw[->, thick]
	(3.8,0.8) -- ++(-1,0.2);

\node at (4.1,0.75) {$i_{\sigma}$};
\node at (2.9,1.2) {$t_{\sigma}$};
\node at (3.4,1.05) {$\sigma$};

\node at (4.2,0) {$1$};
\node at (1,2.2) {$\tau$};

\node[inner sep=0.5, draw, fill=white, shape=circle] at (2.6,-0.2) {\tiny 2};

\end{scope}

\end{tikzpicture}
\caption{Minimal type I tilings on $T_{\tau}$. }
\label{tiling1b}
\end{figure}
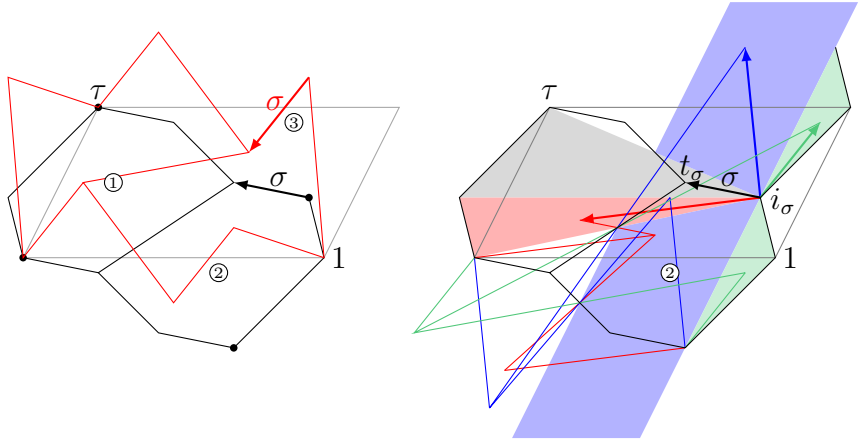

The fundamental domain of the minimal tiling and the four tiles are illustrated in the first of Figure \ref{tiling1b}. The fundamental domain can be any parallelogram. Therefore the torus can be $T_{\tau}$ for any $\tau$ in the upper half plane.

\begin{proof}
The hexagon is the tile $t_1$ in Figure \ref{tiling1type}. We indicate $\bar{0}$ and $\bar{1}$ by green and blue, and indicate $\bar{2}$ and $\bar{5}$ (of equal length) by red. 

The corner $1_1$ of $t_1$ is at a degree $3$ vertex. Let $t_2,t_3$ be the other two tiles at the vertex, on the other side of $\bar{1}_1,\bar{0}_1$. The edge between $t_2,t_3$ is adjacent to $\bar{0}$ and $\bar{1}$. Since $\bar{0}$ and $\bar{1}$ have different lengths from all the others, this implies the edge is red. This determines (all sides and corners of) $t_2,t_3$.

\begin{figure}[htp]
\centering
\begin{tikzpicture}[>=latex,scale=1]

\pgfmathsetmacro{\ra}{(3*sqrt(3)/10)};

\foreach \a in {-1,0,1}
\foreach \c in {1,0,-1}
{
\begin{scope}[xshift=1.8*\a cm, yshift=2*\ra*\c cm]

\draw
	(120:0.6) -- (180:0.6) -- (240:0.6)
	(0:1.2) -- ++(60:0.6) -- ++(120:0.6);

\draw[red]
	(120:0.6) -- (60:0.6) 
	(-60:0.6) -- (240:0.6)
	(60:1.2) -- ++(0:0.6)
	(0:0.6) -- ++(0:0.6);
	
\draw[emerald]
	(-60:0.6) -- (0:0.6)
	(60:0.6) -- (60:1.2);

\draw[blue]
	(60:0.6) -- (0:0.6);

\foreach \b in {0,...,5}
{
\node at (60*\b-60:0.45) {\tiny \b};
\node[shift={(30:2*\ra)}] at (60*\b+120:0.45) {\tiny \b};
}

\end{scope}
}

\node[inner sep=0.5, draw, shape=circle] at (0,0) {\tiny 1};
\node[inner sep=0.5, draw, shape=circle] at (-30:2*\ra) {\tiny 3};
\node[inner sep=0.5, draw, shape=circle] at (30:2*\ra) {\tiny 2};
\node[inner sep=0.5, draw, shape=circle] at (0,2*\ra) {\tiny 4};
\node[inner sep=0.5, draw, shape=circle] at (0,-2*\ra) {\tiny 5};

\end{tikzpicture}
\caption{Tiling by generic hexagon of type I.}
\label{tiling1type}
\end{figure}
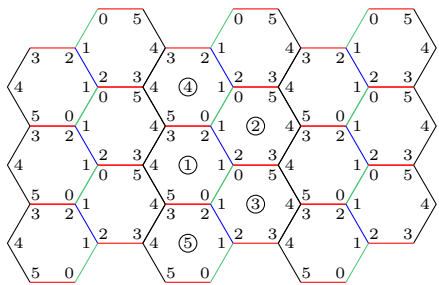

The argument that derives $t_2,t_3$ from $t_1$ can be applied to $t_2,t_3$. Then we determine $t_4,t_5$. More repetitions of the argument determine the two column strip bounded by normal edges in the middle of Figure \ref{tiling1type}. Then the tiling is the union of such strips. The boundary of the strip consists of $\bar{3},\bar{4}$. Since $\bar{3},\bar{4}$ have different lengths, the way to glue the strips is unique, as in Figure \ref{tiling1type}.

Let ${\mc T}_{\text{I}}$ be the unique tiling in Proposition \ref{combo1}. It is schematically given by Figure \ref{tiling1type}, and more realistically given by Figure \ref{tiling1a}. We also use ${\mc T}_{\text{I}}$ to denote the collection of all tiles. Let $G({\mc T}_{\text{I}})$ be the symmetry group of the tiling. We fix one tile $t_1\in {\mc T}_{\text{I}}$ and get a map 
\[
g\mapsto gt_1\colon G({\mc T}_{\text{I}})\to {\mc T}_{\text{I}}.
\]
If $t$ is any tile in ${\mc T}_{\text{I}}$. Then there is the unique isometry $g$ sending $t_1$ to $t$. The uniqueness follows from the matching of the $\bar{0}$ and $\bar{1}$ sides of $t_1$ and $t$. Then the unique tiling implies that $g$ sends the six tiles around $t_1$ to the corresponding six tiles around $t$. Then we further apply the unique tiling to the tiles around each of the six tiles, and find $g$ sending more tiles near $t_1$ to the corresponding tiles near $t$. Repeating the argument shows that $g$ is a symmetry of ${\mc T}_{\text{I}}$.

\begin{figure}[htp]
\centering
\begin{tikzpicture}[>=latex,scale=0.4]

\foreach \a in {-1,...,2}
\foreach \b in {-1,0,1}
\fill[gray!30, shift={(4*\a cm + 1*\b cm,2*\b cm)}]
	(0,0) -- (-0.2,0.8) -- (1,2) -- (2,1.8) -- (2.8,1) -- (1,-0.2);

\foreach \a in {-1,...,2}
\foreach \b in {-1,0,1}
{
\begin{scope}[shift={(4*\a cm + 1*\b cm,2*\b cm)}]	

\draw
	(1.8,-1) -- (1,-0.2) -- (2.8,1) -- (2,1.8);

\draw[red]
	(3.8,0.8) -- ++(-1,0.2)
	(0,0) -- ++(1,-0.2)
	(1,2) -- ++(1,-0.2)
	(2.8,-1.2) -- ++(-1,0.2);

\draw[blue]
	(0,0) -- (-0.2,0.8)
	(4,0) -- (3.8,0.8);

\draw[emerald]
	(-0.2,0.8) -- (1,2)
	(2.8,-1.2) -- (4,0);

\end{scope}
}

\node at (1.9,0.4) {\scriptsize $\lozenge$};

\node at (1.1,1.7) {\tiny $0$};
\node at (0.2,0.8) {\tiny $1$};
\node at (0.3,0.25) {\tiny $2$};
\node at (0.95,0.1) {\tiny $3$};
\node at (2.35,1.05) {\tiny $4$};
\node at (1.9,1.5) {\tiny $5$};

\node at (2.75,-0.9) {\tiny $0$};
\node at (3.6,0) {\tiny $1$};
\node at (3.5,0.6) {\tiny $2$};
\node at (2.85,0.7) {\tiny $3$};
\node at (1.45,-0.2) {\tiny $4$};
\node at (1.9,-0.7) {\tiny $5$};

\node[inner sep=0.5, draw, fill=white, shape=circle] at (1.2,1) {\tiny 1};
\node[inner sep=0.5, draw, fill=white, shape=circle] at (2.6,-0.2) {\tiny 2};
%\node[inner sep=0.5, draw, fill=white, shape=circle] at (3.6,1.8) {\tiny 3};

\end{tikzpicture}
\caption{A real example of the tiling ${\mc T}_{\text{I}}$.}
\label{tiling1a}
\end{figure}
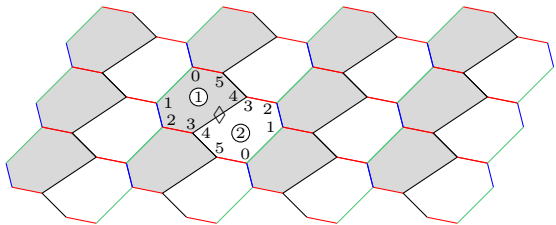

The translations sending gray tiles in Figure \ref{tiling1a} to gray tiles form a lattice subgroup $L\cong{\bb Z}^2$ of $G({\mc T}_{\text{I}})$. The rotation $\rho$ around the middle point $\lozenge$ of $\bar{3}_1=\bar{3}_2$ sends $t_1$ to $t_2$. Then the coset $L\rho$ is all the rotations around the middle points of all edges $\bar{3}$. Under the map $g\mapsto gt_1\colon G({\mc T}_{\text{I}})\to {\mc T}_{\text{I}}$, the image of $L$ is all the gray tiles, and the image of $L\rho$ is all the white tiles. Therefore $G({\mc T}_{\text{I}})$ is transitive on all tiles. Since the generic type I prototile is not symmetric, we conclude the map is a one-to-one correspondence, and we have 
\[
{\mc T}_{\text{I}}
=G({\mc T}_{\text{I}})=L\sqcup L\rho
=L\rtimes\langle\rho\rangle
\cong {\bb Z}^2\rtimes{\bb Z}_2.
\]
The first equality is set-theoretical (and depends on the choice of $t_1$), and the other equalities are group theoretical.

A type I torus tiling is the quotient of ${\mc T}_{\text{I}}$ by a subgroup $G\sub G({\mc T}_{\text{I}})$ that has finite index and acts freely on the plane. Since an orientation preserving isometry of ${\bb C}$ has no fixed point if and only if it is a translation, the lattice subgroup $L$ is exactly all the symmetries of ${\mc T}_{\text{I}}$ without fixed points. Therefore the condition on $G$ means exactly $G\sub L$. Then the type I torus tiling ${\mc T}_{\text{I}}/G$ on ${\bb C}/G$ covers the tiling ${\mc T}_{\text{I}}/L$ on ${\bb C}/L$. This implies ${\mc T}_{\text{I}}/L$, with two tiles, is the minimal type I torus tiling.

Next we identify the torus ${\bb C}/L$ with $T_{\tau}$. In Figure \ref{tiling1a}, we pick two tiles $t_1,t_2$ sharing an edge $\bar{3}$. Let $\alpha$ be the complex number representing the vector from $2_1$ to $1_2$, and let $\beta$ be the complex number representing the vector from $2_1$ to $0_1$. Then the parallelogram spanned by $\alpha,\beta$ is a fundamental domain of $L$. The minimal tiling ${\mc T}_{\text{I}}/L$ is on the torus
\[
{\bb C}/L={\bb C}/({\bb Z}\alpha+{\bb Z}\beta)=T_{\tau},\quad
\tau=\tfrac{\beta}{\alpha}.
\]
Then a minimal type I torus tiling ${\mc T}_{\text{I}}/L$ is naturally associated with a torus $T_{\tau}$.
\end{proof}

We describe the moduli space ${\mc M}^{\text{I}}_{\tau}$ of all minimal type I tilings on the torus $T_{\tau}$. In Figure \ref{tiling1b}, we see such a tiling is determined by the choice of a free vector $\sigma$. The free vector has the initial point $i_{\sigma}$ and the terminal point $t_{\sigma}$. The moduli space ${\mc M}^{\text{I}}_{\tau}$ is the set of all $\sigma$, such that the hexagon $t_2$ constructed from $\sigma$ is simple (i.e., the boundary of the hexagon does not intersect itself).

In the second of Figure \ref{tiling1b}, we fix the initial point $i$ and try to find the possible locations of $t$. Using $i$ as the cone point, we divide the space into several cone regions. For $t$ in various cone regions, we find $t$ should be in the gray, red, blue, and green regions in order for the hexagon to be simple. The union of these regions is the overall possible locations $M^{\rm I}(i)$ of $t$ for the given $i$. Let $I_{\tau}$ be the possible locations of $i$. Then the moduli space
\[
{\mc M}^{\rm I}_{\tau}
=\{i\xrightarrow{\sigma}t\colon i \in I_{\tau},\; t\in M^{\rm I}(i)\}.
\]

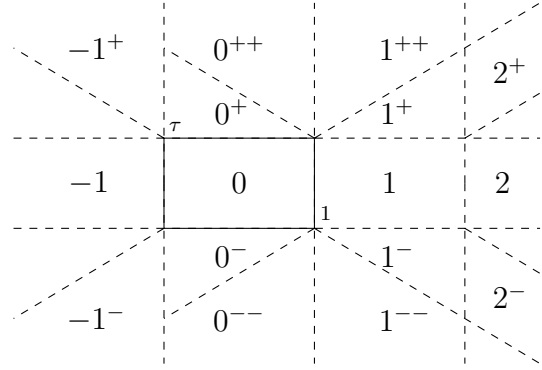
\begin{figure}[htp]
\centering
\begin{tikzpicture}[>=latex,scale=1]

\foreach \a in {1,-1}
\draw[dashed, yscale=\a]
	(-4,0.6) -- (3,0.6)
	(-2,0) -- ++(0,2.4)
	(0,0) -- ++(0,2.4)
	(2,0) -- ++(0,2.4)
	(-2,0.6) -- (-4,1.8)
	(0,0.6) -- (-2,1.8)
	(0,0.6) -- (3,2.4)
	(2,0.6) -- (3,1.2);

\draw
	(0,-0.6) rectangle (-2,0.6);

\node at (0.15,-0.4) {\scriptsize 1};
\node at (-1.85,0.75) {\scriptsize $\tau$};

\node at (-1,0) {0};
\node at (-1.1,0.95) {$0^+$};
\node at (-1,1.8) {$0^{++}$};
\node at (-1.1,-0.95) {$0^-$};
\node at (-1,-1.8) {$0^{--}$};

\node at (-3,0) {$-1$};
\node at (-2.9,1.8) {$-1^+$};
\node at (-2.9,-1.8) {$-1^-$};

\node at (1,0) {1};
\node at (1.1,0.95) {$1^+$};
\node at (1.2,1.8) {$1^{++}$};
\node at (1.1,-0.95) {$1^-$};
\node at (1.2,-1.8) {$1^{--}$};

\node at (2.5,0) {$2$};
\node at (2.6,1.5) {$2^+$};
\node at (2.6,-1.5) {$2^-$};

\end{tikzpicture}
\caption{Regions for the initial $i$. }
\label{tiling1c}
\end{figure}

We remark that the description of the moduli space of possible hexagons in the second of Figure \ref{tiling1b} is actually affine invariant. Therefore we give details of the moduli space only for the rectangular fundamental domain. The shape of $M^{I}(i)$ depends on the locations of $i$. In Figure \ref{tiling1c}, we divide the plane into 22 regions, with the region labeled $0$ being the fundamental domain spanned by $1$ and $\tau$. The 16 labeled regions are the locations of $i$, such that the corresponding region $M^{\rm I}(i)$ for $t$ is not empty. In other words, $I_{\tau}$ is the union of the 16 regions.

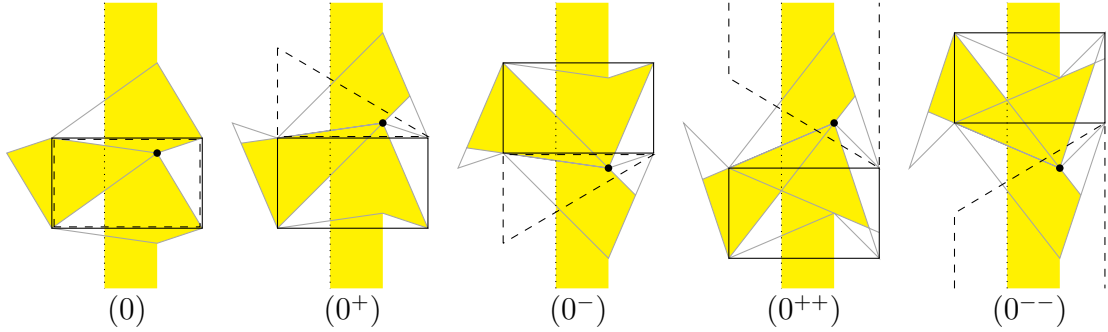
\begin{figure}[htp]
\centering
\begin{tikzpicture}[>=latex,scale=1]

%% 1

\fill[yellow]
	(0.4,0.4) -- (1,0.6) -- (0.4,1.6) -- (0.4,2.4) -- (-0.3,2.4) -- (-0.3,0.5) -- (-1,0.6) -- (-1.6,0.4) -- (-1,-0.6) -- (-0.3,-0.1) -- (-0.3,-1.4) -- (0.4,-1.4) -- (0.4,-0.8) -- (1,-0.6);

\draw[gray!70]
	(-1,-0.6) -- (0.4,-0.8) -- (1,-0.6) -- (0.4,0.4) -- (1,0.6) -- (0.4,1.6) -- (-1,0.6) -- (-1.6,0.4) -- (-1,-0.6) -- (0.4,0.4) -- (-1,0.6);

\draw
	(-1,-0.6) rectangle (1,0.6);

\draw[dashed, scale=0.97]
	(-1,-0.6) rectangle (1,0.6);

\draw[dotted]
	(-0.3,2.4) -- (-0.3,-1.4);
	
\fill (0.4,0.4) circle (0.05);

\node at (0,-1.7) {(0)};
\node at (3,-1.7) {($0^+$)};
\node at (6,-1.7) {($0^-$)};
\node at (9,-1.7) {($0^{++}$)};
\node at (12,-1.7) {($0^{--})$};

%% 2

\foreach \a in {1,-1}
{
\begin{scope}[shift={(4.5cm -1.5*\a cm, 0.5 cm-0.5 *\a cm)}, yscale=\a]

\fill[yellow]
	(0.4,0.8) -- (0.76,1.16) -- (0.4,2) -- (0.4,2.4) -- (-0.3,2.4) -- (-0.3,0.7) -- (-1,0.6) -- (-1.485,0.531) -- (-1,-0.6) -- (-0.3,0.1) -- (-0.3,-1.4) -- (0.4,-1.4) -- (0.4,-0.4) -- (1,-0.6);

\draw[gray!70]
	(-1,-0.6) -- (0.4,-0.4) -- (1,-0.6) -- (0.4,0.8) -- (1,0.6) -- (0.4,2) -- (-1,0.6) -- (-1.6,0.8) -- (-1,-0.6) -- (0.4,0.8) -- (-1,0.6) 
	(-1.485,0.531) -- (0.4,0.8) -- (0.76,1.16);

\draw
	(-1,-0.6) rectangle (1,0.6);

\draw[dashed]
	(-1,1.8) -- (-1,0.62) -- (1,0.62) -- (-1,1.8);

\draw[dotted]
	(-0.3,2.4) -- (-0.3,-1.4);
		
\fill (0.4,0.8) circle (0.05);

\end{scope}
}

%% 3
\foreach \a in {1,-1}
{
\begin{scope}[shift={(10.5cm -1.5*\a cm, 0.5 cm-0.9 *\a cm)}, yscale=\a]

\fill[yellow]
	(0.4,1.2) -- (0.68,1.56) -- (0.4,2.4) -- (0.4,2.8) -- (-0.3,2.8) -- (-0.3,0.9) -- (-1.35,0.45) -- (-1,-0.6) -- (-0.3,0.3) -- (-0.3,-1) -- (0.4,-1) -- (0.4,0) -- (0.867,-0.2);

\draw[gray!70]
	(-1,-0.6) -- (0.4,0) -- (1,-0.6) -- (0.4,1.2) -- (1,0.6) -- (0.4,2.4) -- (-1,0.6) -- (-1.6,1.2) -- (-1,-0.6) -- (0.4,1.2) -- (-1,0.6)
	(-1.35,0.45) -- (0.4,1.2) -- (0.68,1.56)
	(-1,0.6) -- (1,-0.257);

\draw
	(-1,-0.6) rectangle (1,0.6);

\draw[dashed]
	(-1,2.8) -- (-1,1.8) -- (1,0.6) -- (1,2.8);

\draw[dotted]
	(-0.3,2.8) -- (-0.3,-1);
		
\fill (0.4,1.2) circle (0.05);

\end{scope}
}

\end{tikzpicture}
\caption{Region $M^{\rm I}(i)$ of possible $t$, for $i=\bullet$ in regions labeled by $0^*$. }
\label{tiling1_moduli1}
\end{figure}

In Figure \ref{tiling1_moduli1}, the yellow regions show $M^{\rm I}(i)$ for the case $i$ is in the regions labeled by $0^*$. Part of the boundary of $M^{\rm I}(i)$ is the dotted middle line between $i$ (indicated by $\bullet$) and the left side $\tau$ of the fundamental domain. 

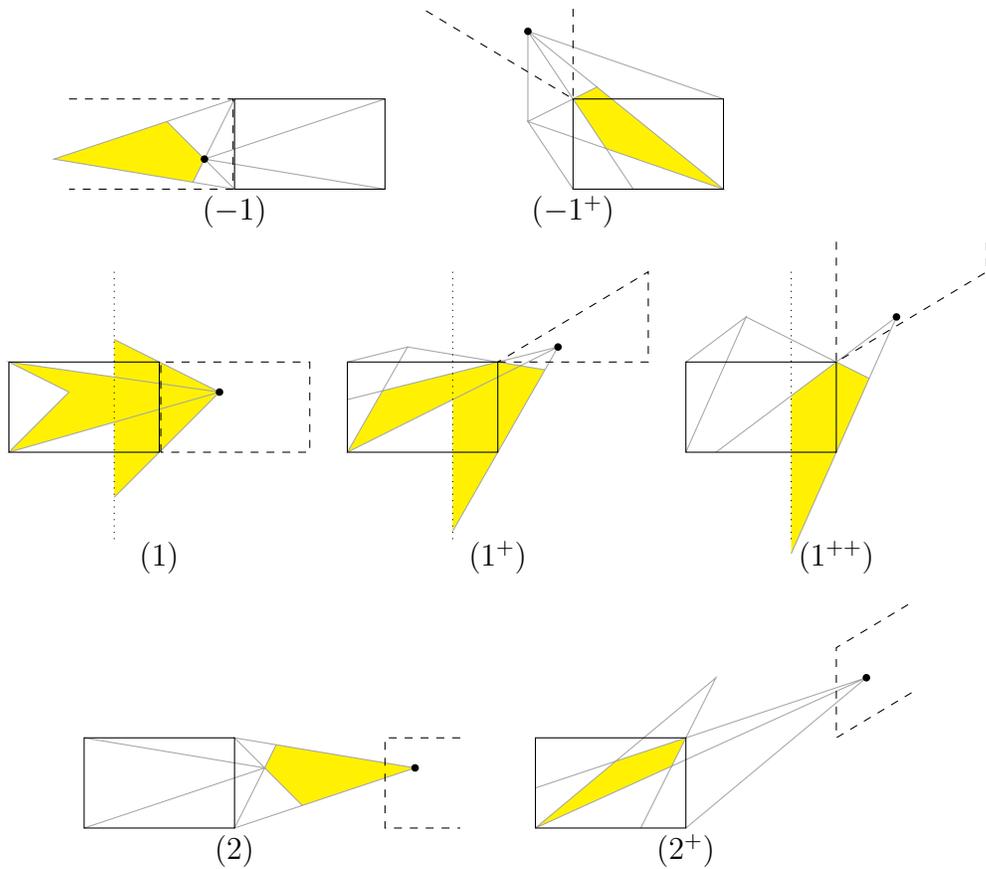
\begin{figure}[htp]
\centering
\begin{tikzpicture}[>=latex,scale=1]

%% -1

\begin{scope}[shift={(3cm, 3.5cm)}]

\fill[yellow]
	(-3.4,-0.2) -- (-1.554,-0.508) -- (-1.4,-0.2) -- (-1.9,0.3);

\draw[gray!70]
	(1,-0.6) -- (-1.4,-0.2) -- (1,0.6)
	(-1.9,0.3) -- (-1,-0.6) -- (-3.4,-0.2) -- (-1,0.6) -- (-1.554,-0.508);

\draw
	(-1,-0.6) rectangle (1,0.6);

\draw[dashed]
	(-3.2,0.6) -- (-1.02,0.6) -- (-1.02,-0.6) -- (-3.2,-0.6);
	
\fill (-1.4,-0.2) circle (0.05);

\node at (-1,-0.9) {($-1$)};

\begin{scope}[xshift=4.5cm]

\fill[yellow]
	(-1,0.6) -- (-0.682,0.759) -- (1,-0.6) -- (-0.56,-0.06);

\draw[gray!70]
	(-1,-0.6) -- (-1.6,0.3) -- (1,-0.6) -- (-1.6,1.5) -- (-0.2,-0.6)
	(-0.682,0.759) -- (-1.6,0.3) -- (-1.6,1.5) -- (1,0.6);

\draw
	(-1,-0.6) rectangle (1,0.6);

\draw[dashed]
	(-1,1.8) -- (-1,0.6) -- (-3,1.8) -- (-3,1.8);
	
\fill (-1.6,1.5) circle (0.05);

\node at (-1,-0.9) {($-1^+$)};

\end{scope}

\end{scope}

%% 1

\fill[yellow]
	(-0.2,0.2) -- (-1,0.6) -- (0.4,0.4) -- (0.4,0.9) -- (1.8,0.2)-- (0.4,-1.2) -- (0.4,-0.2) --  (-1,-0.6) -- (-0.2,0.2);

\draw[gray!70]
	(0.4,-1.2) -- (1.8,0.2)--  (0.4,0.9)
	(-0.2,0.2) -- (-1,-0.6) -- (1.8,0.2)--  (-1,0.6) -- (-0.2,0.2);

\draw
	(-1,-0.6) rectangle (1,0.6);

\draw[dashed]
	(1.02,-0.6) rectangle (3,0.6);

\draw[dotted]
	(0.4,1.8) -- (0.4,-1.8);
	
\fill (1.8,0.2) circle (0.05);

\node at (1,-2) {(1)};

\begin{scope}[xshift=4.5cm]

\fill[yellow]
	(-1,-0.6) -- (-0.533,0.217) -- (1,0.6) -- (1.626,0.496) -- (0.4,-1.65) -- (0.4,0.1);

\draw[gray!70]
	(-1,0.6) -- (-0.2,0.8) --  (-1,-0.6) -- (1.8,0.8)
	(-1,0.1) -- (1.8,0.8) --  (0.4,-1.65)
	(-0.2,0.8) -- (1.626,0.496);

\draw
	(-1,-0.6) rectangle (1,0.6);

\draw[dashed]
	 (1,0.6) -- (3,0.6) -- (3,1.8) -- (1,0.6);

\draw[dotted]
	(0.4,1.8) -- (0.4,-1.8);
	
\fill (1.8,0.8) circle (0.05);

\node at (1,-2) {($1^{+}$)};

\end{scope}

\begin{scope}[xshift=9cm]

\fill[yellow]
	(1,0.6) -- (1.44,0.38) -- (0.4,-1.95) -- (0.4,0.15);

\draw[gray!70]
	(-1,0.6)-- (-0.2,1.2) -- (-1,-0.6)
	(-0.6,-0.6) -- (1.8,1.2) --  (0.4,-1.95)
	(-0.2,1.2) -- (1.44,0.38);

\draw
	(-1,-0.6) rectangle (1,0.6);

\draw[dashed]
	(1,2.2) -- (1,0.6) -- (3,1.8) -- (3,2.2);

\draw[dotted]
	(0.4,1.8) -- (0.4,-1.8);
	
\fill (1.8,1.2) circle (0.05);

\node at (1,-2) {($1^{++}$)};

\end{scope}

%% 2

\begin{scope}[shift={(1cm, -5cm)}]

\fill[yellow]
	(3.4,0.2) -- (1.9,-0.3) -- (1.4,0.2) -- (1.554,0.508);

\draw[gray!70]
	(-1,0.6) -- (1.4,0.2) --  (-1,-0.6)
	(1.9,-0.3) -- (1,0.6) -- (3.4,0.2) --  (1,-0.6) -- (1.554,0.508)
	;

\draw
	(-1,-0.6) rectangle (1,0.6);

\draw[dashed]
	(4,0.6) -- (3,0.6) -- (3,-0.6) -- (4,-0.6);
	
\fill (3.4,0.2) circle (0.05);

\node at (1,-0.9) {(2)};

\begin{scope}[xshift=6cm]

\fill[yellow]
	(1,0.6) -- (0.812,0.224) -- (-1,-0.6) -- (0.067,0.289);

\draw[gray!70]
	(-1,-0.6) -- (1.4,1.4) --  (0.4,-0.6)
	(-1,-0.6) -- (3.4,1.4) --  (1,-0.6)
	(3.4,1.4) -- (-1,-0.067);

\draw
	(-1,-0.6) rectangle (1,0.6);

\draw[dashed]
	(4,1.2) -- (3,0.6) -- (3,1.8) -- (4,2.4);
	
\fill (3.4,1.4) circle (0.05);

\node at (1,-0.9) {($2^+$)};

\end{scope}

\end{scope}

\end{tikzpicture}
\caption{Region $M^{\rm I}(i)$ of possible $t$, for $i=\bullet$ in regions labeled by $-1^*$, $1^*$, $2^*$. }
\label{tiling1_moduli2}
\end{figure}

We observe that $M^{\rm I}(i)$ for $0^+$ and $0^-$ are related by vertical flip. Specifically, suppose $i$ in the region labeled $0^+$ and $i'$ in the region labeled $0^-$ are related by the vertical flip that preserves the fundamental rectangular domain, then $M^{\rm I}(i)$ and $M^{\rm I}(i')$ are also related by the same vertical flip. The same observation holds for $0^{++}$ and $0^{--}$. In fact, this also holds if $0$ is replaced by $-1$, $1$ or $2$. Therefore in Figure \ref{tiling1_moduli2}, we give $M^{\rm I}(i)$ only for $i$ in regions labeled $k,k^+,k^{++}$.

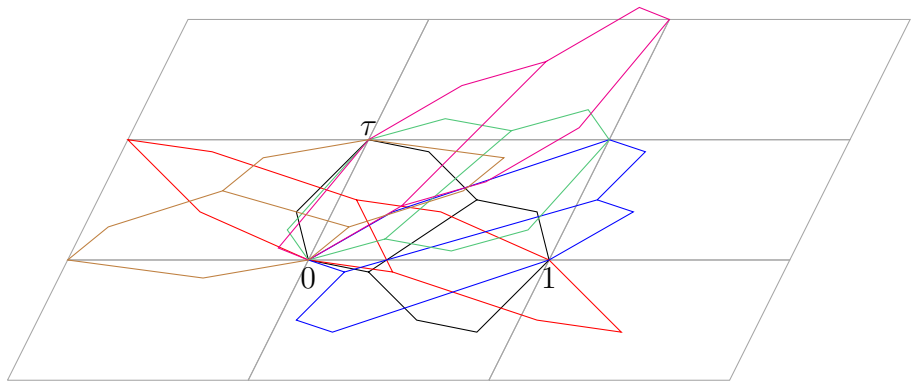
\begin{figure}[htp]
\centering
\begin{tikzpicture}[>=latex,scale=0.8]

\foreach \a in {-1,0,1}
\foreach \b in {-1,0,1}
\draw[gray!70, shift={(4*\a cm+1*\b cm, 2*\b cm)}]
	(0,0) -- (4,0) -- (5,2) -- (1,2) -- (0,0); 	

\draw
	(1,-0.2) -- (0,0) -- (-0.2,0.8) -- (1,2) -- (2,1.8) -- (2.8,1) -- (1,-0.2) -- (1.8,-1) -- (2.8,-1.2) -- (4,0) -- (3.8,0.8) -- (2.8,1);
	
\draw[red, y={(-2cm,1cm)}]
	(1,-0.2) -- (0,0) -- (-0.2,0.8) -- (1,2) -- (2,1.8) -- (2.8,1) -- (1,-0.2) -- (1.8,-1) -- (2.8,-1.2) -- (4,0) -- (3.8,0.8) -- (2.8,1);

\draw[blue, y={(2cm,1cm)}]
	(1,-0.2) -- (0,0) -- (-0.2,0.8) -- (1,2) -- (2,1.8) -- (2.8,1) -- (1,-0.2) -- (1.8,-1) -- (2.8,-1.2) -- (4,0) -- (3.8,0.8) -- (2.8,1);

\draw[emerald, x={(1.25cm,0.5cm)}, y={(-0.125cm,0.75cm)}]
	(1,-0.2) -- (0,0) -- (-0.2,0.8) -- (1,2) -- (2,1.8) -- (2.8,1) -- (1,-0.2) -- (1.8,-1) -- (2.8,-1.2) -- (4,0) -- (3.8,0.8) -- (2.8,1);

\draw[magenta, x={(1.5cm,1cm)}, y={(-0.25cm,0.5cm)}]
	(1,-0.2) -- (0,0) -- (-0.2,0.8) -- (1,2) -- (2,1.8) -- (2.8,1) -- (1,-0.2) -- (1.8,-1) -- (2.8,-1.2) -- (4,0) -- (3.8,0.8) -- (2.8,1);

\draw[brown, x={(0.25cm,0.5cm)}, y={(-2.125cm,-0.25cm)}]
	(1,-0.2) -- (0,0) -- (-0.2,0.8) -- (1,2) -- (2,1.8) -- (2.8,1) -- (1,-0.2) -- (1.8,-1) -- (2.8,-1.2) -- (4,0) -- (3.8,0.8) -- (2.8,1);
			
%\draw[thick] (3.8,0.8) -- ++(-1,0.2);

\node at (4,-0.3) {$1$};
\node at (1,2.2) {$\tau$};
\node at (0,-0.3) {$0$};
%\node at (3.4,1.05) {$\sigma$};

\end{tikzpicture}
\caption{$\mu{\mc T}$ for $\mu=
{\color{red} \begin{pmatrix} 1 & -1 \\ 0 & 1 \end{pmatrix} }, 
{\color{blue} \begin{pmatrix} 1 & 1 \\ 0 & 1 \end{pmatrix} }, 
{\color{emerald} \begin{pmatrix} 1 & 0 \\ 1 & 1 \end{pmatrix} }, 
{\color{magenta} \begin{pmatrix} 1 & 0 \\ 2 & 1 \end{pmatrix} }, 
{\color{brown} \begin{pmatrix} 0 & -1 \\ 1 & 0 \end{pmatrix} }$.}
\label{tiling1d}
\end{figure}

Next we fix an isometric class of the torus, and consider all minimal type I tilings for the class. Let $T_{\tau}={\bb C}/\Lambda_{\tau}$ be a representative in the class. Then minimal tilings on $T_{\tau}$ also include tilings in ${\mc M}^{\text{I}}_{\tau'}$, where $\tau'=\mu(\tau)$ for some $\mu\in SL(2,{\bb Z})$. This means the moduli space of all minimal tilings on $T_{\tau}$ is 
\[
\cup_{\tau'\sim\tau}{\mc M}^{\text{I}}_{\tau'}
=SL(2,{\bb Z}){\mc M}^{\text{I}}_{\tau}
=\{\mu{\mc T}\colon {\mc T}\in {\mc M}^{\text{I}}_{\tau},\; \mu\in SL(2,{\bb Z})\}.
\]
The action $\mu{\mc T}$ literally applies the linear transformation $\mu$ to ${\mc T}$, with respect to the basis $1,\tau$. In Figure \ref{tiling1d}, the black tiles form a minimal tiling ${\mc T}$ of $T_{\tau}$, and colored tilings are transformed tilings $\mu{\mc T}$ on the same $T_{\tau}$ for various $\mu\in SL(2,{\bb Z})$. Since the transformed tiling $\mu{\mc T}$ is characterised by a new basis of the lattice $\Lambda_{\tau}$, we see that $SL(2,{\bb Z}){\mc M}^{\text{I}}_{\tau}=\sqcup_{\mu\in SL(2,{\bb Z})}\mu{\mc M}^{\text{I}}_{\tau}$ is topologically a disjoint union.

\subsection{Generic Type II Tiling}

\begin{theorem}\label{combo2}
A type ${\rm II}$ hexagonal prototile satisfies
\[
[0]+[1]+[3]=[2]+[4]+[5]=2\pi,\quad
|\bar{1}|=|\bar{3}|,\quad
|\bar{2}|=|\bar{5}|.
\]
Suppose it further satisfies the following generic conditions:
\begin{itemize}
\item Each of $\bar{0},\bar{4}$ has different lengths from the other five sides.
\item $|\bar{1}|\ne |\bar{2}|$.
\item $[2]\ne[3]$.
\end{itemize}
Then its tiling of the plane, such that all vertices have degree $3$, is uniquely given by Figures \ref{tiling2type} and \ref{tiling2a}. The tiling is isohedral, periodic, and covers a minimal torus tiling with four tiles. 
\end{theorem}

The fundamental domain of the minimal tiling and the four tiles are illustrated in Figure \ref{tiling2b}. The fundamental domain can be any rectangle. Therefore the torus can be $T_{\tau}$ for any purely imaginary $\tau$ in the upper half plane.

\begin{proof}
The hexagon is the tile $t_1$ in Figure \ref{tiling2type}. We indicate the sides $\bar{1}$ and $\bar{3}$ (of equal length) by red, and indicate the sides $\bar{2}$ and $\bar{5}$ (of equal length) by blue. 

Let $t_2$ be the tile on the other side of $\bar{0}_1$. Since $\bar{0}$ has different length from the other five, we get $\bar{0}_1=\bar{0}_2$. If the direction of $\bar{0}_2$ is not as indicated, then we get two blue sides at the degree $3$ vertex $0_11_2\cdots$, contradicting the non-adjacency of blue sides in the hexagon. Therefore the direction of $\bar{0}_2$ is as indicated. This determines $t_2$. By the same reason, we determine $t_3$. Then we get $t_4$, in which we already know two red sides and the blue side between them. This implies $\bar{5}_1=\bar{2}_4$. If the direction of $\bar{2}_4$ is not as indicated, then the vertex between $t_1,t_2,t_4$ is $0_11_22_4$. By the angle sum of the vertex and $[0]+[1]+[3]=2\pi$, we get $[2]=[3]$, a contradiction. Therefore the direction of $\bar{2}_4$ is not as indicated. This determines $t_4$. 

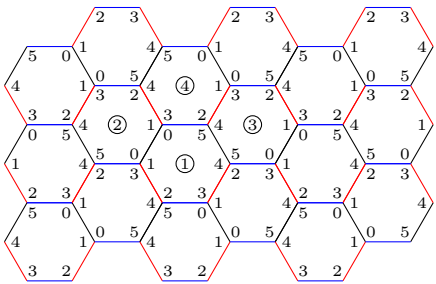
\begin{figure}[htp]
\centering
\begin{tikzpicture}[>=latex,scale=1]

\pgfmathsetmacro{\ra}{(3*sqrt(3)/10)};

\foreach \a in {-1,0,1}
\foreach \c in {1,0,-1}
{
\begin{scope}[xshift=1.8*\a cm, yshift=2*\ra*\c cm]

\draw
	(0:0.6) -- (60:0.6)
	(180:0.6) -- (120:0.6)
	(0:1.2) -- ++(60:0.6);

\draw[red]
	(0:0.6) -- (-60:0.6)
	(180:0.6) -- (240:0.6)
	(60:0.6) -- (60:1.2)
	(1.2,2*\ra) -- ++(-60:0.6);
		
\draw[blue]
	(120:0.6) -- (60:0.6)
	(240:0.6) -- (-60:0.6)
	(0:0.6) -- (0:1.2)
	(60:1.2) -- ++(0:0.6);

\end{scope}
}

\foreach \a in {-1,0,1}
\foreach \b in {0,...,5}
{
\begin{scope}[xshift=1.8*\a cm]

\node at (60*\b+120:0.45) {\tiny \b};
\node[shift={(30:2*\ra)}] at (60*\b-60:0.45) {\tiny \b};

\foreach \c in {1,-1}
{
\begin{scope}[yshift=2*\ra*\c cm]

\node at (-60*\b+60:0.45) {\tiny \b};
\node[shift={(30:2*\ra)}] at (-60*\b-120:0.45) {\tiny \b};

\end{scope}
}

\end{scope}
}

\node[inner sep=0.5, draw, shape=circle] at (0,0) {\tiny 1};
\node[inner sep=0.5, draw, shape=circle] at (150:2*\ra) {\tiny 2};
\node[inner sep=0.5, draw, shape=circle] at (30:2*\ra) {\tiny 3};
\node[inner sep=0.5, draw, shape=circle] at (90:2*\ra) {\tiny 4};

\end{tikzpicture}
\caption{Tiling by generic hexagon of type II.}
\label{tiling2type}
\end{figure}

We see that a tile $t_1$ determines tiles $t_2,t_3,t_4$. If we apply the same argument to $t_2,t_3,t_4$, then more tiles are determined. More repetitions give the whole tiling in Figure \ref{tiling2type}.

Let ${\mc T}_{\text{II}}$ be the unique tiling in Proposition \ref{combo2}. It is schematically given by Figure \ref{tiling2type}, and more realistically given by Figure \ref{tiling2a}. The translations among the gray tiles form a subgroup $L\cong {\bb Z}^2$ of the symmetry group $G({\mc T}_{\text{II}})$ of the tiling. Moreover, the rotation $\rho$ with respect to $\lozenge$ by $180^{\circ}$ takes $t_1$ to $t_2$, and the glide reflection $\gamma$ along the green line segment takes $t_1$ to $t_4$, and $\rho\gamma$ takes $t_1$ to $t_3$. Then $L\sqcup L\rho\sqcup L\gamma\sqcup L\rho\gamma$ is transitive on all the tiles. Since the generic type II prototile is not symmetric, by the same argument for ${\mc T}_{\text{I}}$, we conclude these are all the symmetries
\[
G({\mc T}_{\text{II}})=L\sqcup L\rho\sqcup L\gamma\sqcup L\rho\gamma,
\]
and $g\mapsto gt_1\colon G({\mc T}_{\text{II}})\to{\mc T}_{\text{II}}$ is a one-to-one correspondence.

\begin{figure}[htp]
\centering
\begin{tikzpicture}[>=latex,scale=(1/3)]

\foreach \x in {-1,...,2}
\foreach \y in {-1,...,1}
\fill[gray!30, shift={(6*\x cm, 4*\y cm)}]
	(-0.4,1.6) -- (-2,1.2) -- (-3.4,0.4) -- (-2.6,-0.4) -- (-1,-0.8) -- (1,0.8);
	
\node at (0,0) {\scriptsize $\lozenge$};

\foreach \a in {-1,1}
\foreach \x in {-1,...,2}
\foreach \y in {-1,...,1}
{
\begin{scope}[shift={(6*\x cm, 4*\y cm)}, scale=\a]

\draw
	(1,0.8) -- (-1,-0.8)
	(2.6,0.4) -- (3.4,-0.4)
	(-0.4,1.6) -- (0.4,2.4)
	(2,2.8) -- (4,1.2);

\draw[red]
	(1,0.8) -- (-0.4,1.6)
	(2,-1.2) -- (3.4,-0.4)
	(4,1.2) -- (2.6,0.4);

\draw[blue]
	(1,0.8) -- (2.6,0.4)
	(0.4,2.4) -- (2,2.8)
	(0.4,-1.6) -- (2,-1.2);

\end{scope}
}

\foreach \a in {-1,1}
{
\begin{scope}[scale=\a]

\node at (1.1,0.4) {\tiny 0};
\node at (-0.4,-0.8) {\tiny 1};
\node at (0.4,-1.2) {\tiny 2};
\node at (1.9,-0.85) {\tiny 3};
\node at (2.9,-0.3) {\tiny 4};
\node at (2.4,0.1) {\tiny 5};

\end{scope}
}

\node[inner sep=0.5, draw, fill=white, shape=circle] at (-1.3,0.5) {\tiny 1};
\node[inner sep=0.5, draw, fill=white, shape=circle] at (1.3,-0.5) {\tiny 2};
\node[inner sep=0.5, draw, fill=white, shape=circle] at (1.7,1.5) {\tiny 3};
\node[inner sep=0.5, draw, fill=white, shape=circle] at (-1.7,-1.5) {\tiny 4};
\node[inner sep=0.5, draw, fill=white, shape=circle] at (-1.7,2.5) {\tiny 5};
\node[inner sep=0.5, draw, fill=white, shape=circle] at (1.7,-2.5) {\tiny 6};

\draw[thick, emerald]
	(-1.5,0.5) -- (-1.5,-1.5);
		
\end{tikzpicture}
\caption{A real example of the tiling ${\mc T}_{\text{II}}$.}
\label{tiling2a}
\end{figure}
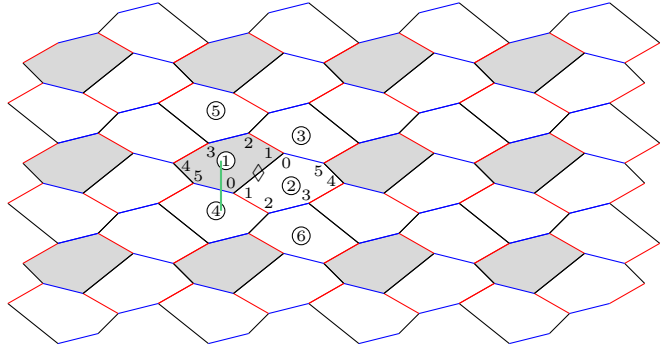

We have $\gamma\rho\gamma=\rho$, and $L\sqcup L\gamma$ is the subgroup of symmetries without fixed points, and
\[
G({\mc T}_{\text{II}})=(L\sqcup L\gamma)\rtimes \langle\rho\rangle.
\]
Then a type II torus tiling is the quotient of ${\mc T}_{\text{II}}$ by a subgroup $G\sub L\sqcup L\gamma$ of finite index. Moreover, $G$ should contain no orientation reversing symmetries because the quotient surface ${\bb C}/G$ would otherwise be the Klein bottle. Therefore $G\sub L$. This implies ${\mc T}_{\text{II}}/L$ on ${\bb C}/L$, with four tiles, is the minimal type II tiling. 

If we include orientation reversing symmetries in $G$, then $G=G_0\sqcup G_0\gamma$, where $G_0=G\cap L$ is all orientation preserving symmetries in $G$. Then we get a type II tiling ${\mc T}_{\text{II}}/G$ on the Klein bottle ${\bb C}/G$. The minimal type II tiling on Klein bottle is then ${\mc T}_{\text{II}}/(L\sqcup L\gamma)$ on ${\bb C}/(L\sqcup L\gamma)$, and has two tiles. 

In the tiling ${\mc T}_{\text{II}}$ in Figure \ref{tiling2a}, we pick an edge $\bar{0}$ shared by two tiles $t_1,t_2$. At the two ends of the edge, we also have tiles $t_3,t_4$ in addition to $t_1,t_2$. Then we get $t_5,t_6$ respectively sharing degree $3$ vertices with $t_1,t_3$ and $t_2,t_4$. Then we connect the middle points of $\bar{0}_3,\bar{0}_4,\bar{0}_5,\bar{0}_6$ to get a rectangle. This rectangle is a fundamental domain of $L$. Therefore the torus is $T_{\tau}={\bb C}/\Lambda_{\tau}$ for a purely imaginary $\tau$.

\begin{figure}[htp]
\centering
\begin{tikzpicture}[>=latex,scale=0.6]

\draw[gray!70]
	(-3,-2) rectangle (3,2);

\draw[gray!70,dashed]	
	(-3,0) -- (3,0)
	(0,-2) -- (0,2);

\foreach \a in {-1,1}
{
\begin{scope}[scale=\a]

\draw
	(1,0.8) -- (-1,-0.8)
	(2.6,0.4) -- (3.4,-0.4)
	(-0.4,1.6) -- (0.4,2.4)
	(2,2.8) -- (4,1.2);

\draw[dashed]
	(-0.4,-2.4) -- (1,-3.2) -- (2.6,-3.6) -- (4,-2.8) -- (2,-1.2);

\draw[red]
	(1,0.8) -- (-0.4,1.6)
	(2,-1.2) -- (3.4,-0.4)
	(4,1.2) -- (2.6,0.4);

\draw[blue, thick]
	(1,0.8) -- (2.6,0.4)
	(0.4,2.4) -- (2,2.8)
	(0.4,-1.6) -- (2,-1.2);

\fill[blue] 
	(2.6,0.4) circle (0.1)
	(0.4,2.4) circle (0.1)
	(0.4,-1.6) circle (0.1);
	
\end{scope}
}

\draw[blue, thick, <-]
	(1,0.8) -- (2.6,0.4);

\node at (3.4,0.2) {$\frac{1}{2}$};
\node at (-3.55,-0.2) {$-\frac{1}{2}$};
\node at (-0.2,2.3) {$\frac{\tau}{2}$};
\node at (0.3,-2.3) {$-\frac{\tau}{2}$};

\node[blue] at (1.9,0.85) {$\sigma$};

\node[inner sep=0.5, draw, fill=white, shape=circle] at (-1.3,0.5) {\tiny 1};
\node[inner sep=0.5, draw, fill=white, shape=circle] at (-1.7,-1.5) {\tiny 4};
\node[inner sep=0.5, draw, fill=white, shape=circle] at (1.3,-0.5) {\tiny 2};
\node[inner sep=0.5, draw, fill=white, shape=circle] at (1.7,1.5) {\tiny 3};
\node[inner sep=0.5, draw, fill=white, shape=circle] at (-1.7,2.5) {\tiny 5};
\node[inner sep=0.5, draw, fill=white, shape=circle] at (1.7,-2.5) {\tiny 6};
		
\end{tikzpicture}
\caption{Moduli space ${\mc M}^{\text{II}}_{\tau}$ of minimal  type II tilings on $T_{\tau}$. }
\label{tiling2b}
\end{figure}
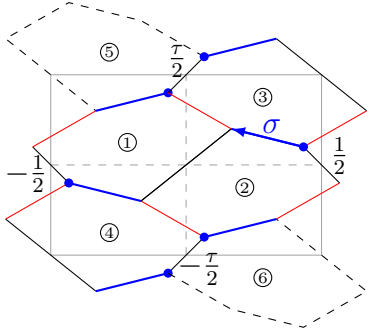

\end{proof}

We choose the fundamental domain of $T_{\tau}$ to be $[-\frac{1}{2},\frac{1}{2}]\times [-\frac{\tau}{2i},\frac{\tau}{2i}]$. A minimal type II tiling on $T_{\tau}$ is given by a free vector $\sigma$. We glide reflect $\sigma$ along the line $x=\frac{1}{4}$ by distance $\pm\frac{\tau}{2i}$ to get two other vectors. Then we rotate the three vectors by $180^{\circ}$ to get three more vectors. More shiftings by multiples of $1$ or $\tau$ give more vectors. Then connecting the ends of these vectors gives tiles $t_1,t_2,t_3,t_4$. The four tiles form a minimal type II tiling of $T_{\tau}$. The moduli space ${\mc M}^{\text{II}}_{\tau}$ is then the space of free vectors $\sigma$, such that the tile $t_2$ constructed above is simple. 

We may still try to describe ${\mc M}^{\text{II}}_{\tau}$ similar to type I tiling. In other words, we fix the initial $i$ of $\sigma$, and describe possible the locations $M^{\rm II}(i)$ of the terminal $t$ of $\sigma$, such that the tile $t_2$ constructed above is simple. The description turns out to be very complicated. We draw the picture by using GeoGebra, and we show two samples here.

\begin{figure}[htp]
\centering
\begin{tikzpicture}[>=latex,scale=1]

\pgftext{
	\includegraphics[scale=0.25]{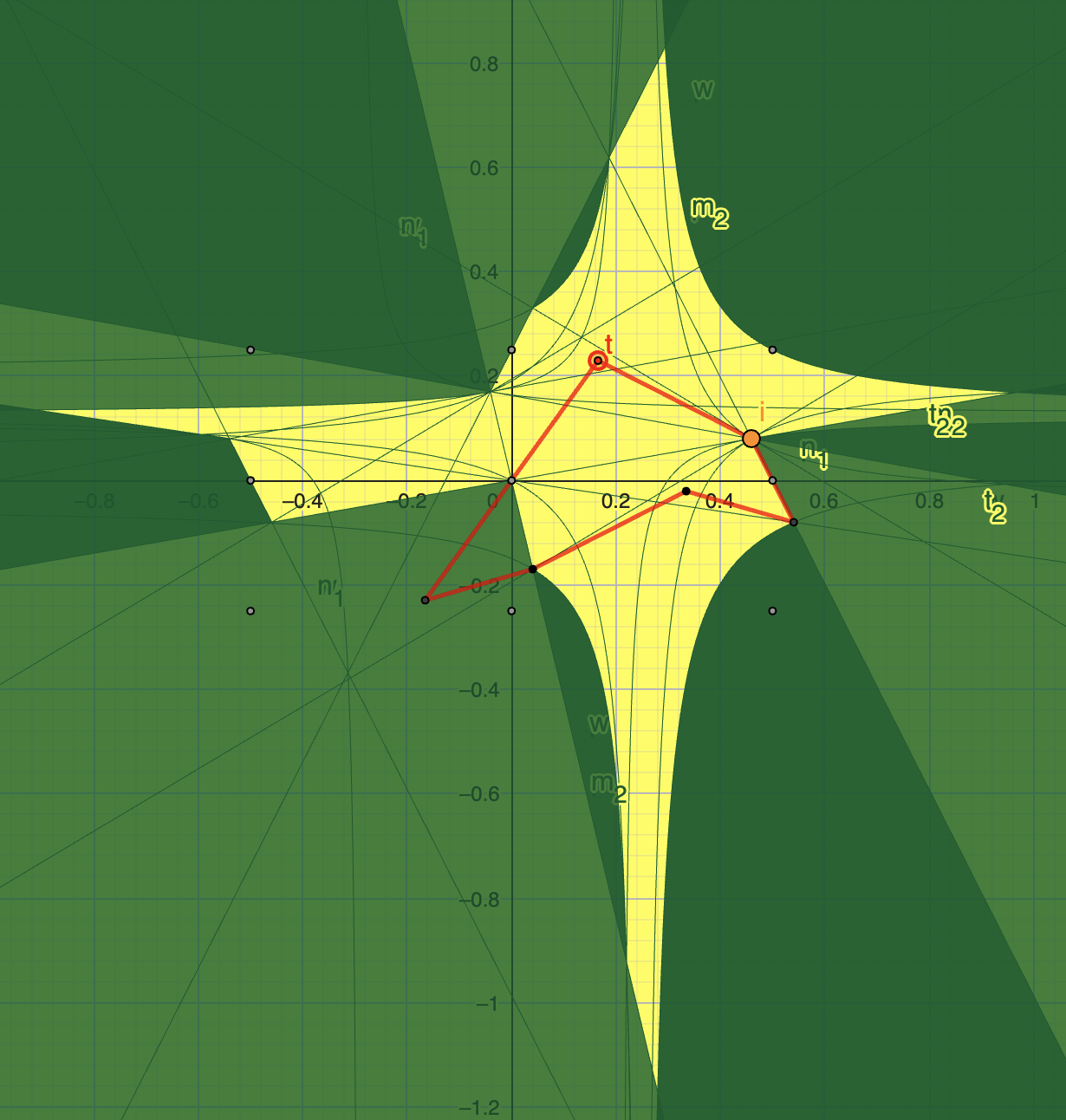}};

\begin{scope}[xshift=7cm]

\pgftext{
	\includegraphics[scale=0.25]{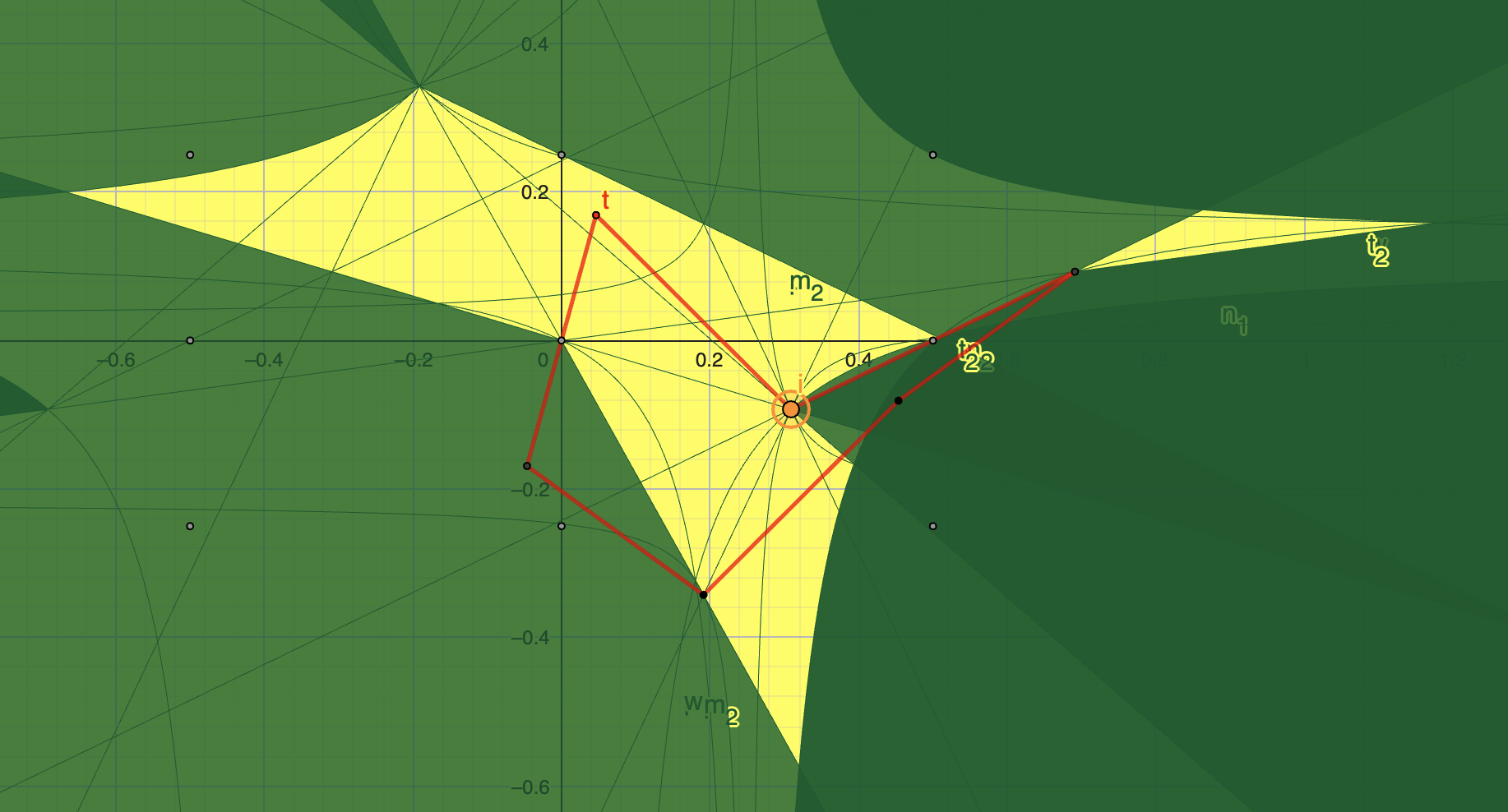}};

\end{scope}

\end{tikzpicture}
\caption{Examples of $M^{\rm II}(i)$, for two locations of $i$. }
\label{tiling2m}
\end{figure}

On the left of Figure \ref{tiling2m}, the initial (orange dot) $i$ lies in the first quadrangle of the fundamental domain. The boundary of the corresponding $M^{\rm II}(i)$ consists of 9 straight line segments and 5 hyperbolic segments. As $i$ moves in the first quadrant, the shape is largely the same, but some segments may disappear. 

On the right of Figure \ref{tiling2m}, the initial $i$ lies in the upper right half of the fourth quadrangle of the fundamental domain. In this case, $M^{\rm II}(i)$ actually has two connected components! 

In general, our experiment with GeoGebra shows that $M^{\rm II}(i)$ is either empty, connected, or has two connected components. The connected components do not have holes. Moreover, the boundary always consists of straight lines and hyperbolic arcs.

Finally, the fundamental domain for type II tiling is always rectangular. Therefore the only way we get different minimal tiling on the same torus like in Figure \ref{tiling1d} is the rotation by $90^{\circ}$ in case $\tau=\mathrm{i}$.

\subsection{Generic Type III Tiling}

\begin{theorem}\label{combo3}
A type ${\rm III}$ hexagonal prototile satisfies 
\[
[1]=[3]=[5]=\tfrac{2}{3}\pi,\quad
|\bar{0}|=|\bar{1}|,\quad
|\bar{2}|=|\bar{3}|,\quad
|\bar{4}|=|\bar{5}|.
\]
Suppose it further satisfies the following generic conditions:
\begin{itemize}
\item $|\bar{0}|\ne |\bar{2}|\ne |\bar{4}|$.
\item $[2],[4],[6]\ne\frac{2}{3}\pi$.
\end{itemize}
Then its tiling of the plane, such that all vertices have degree $3$, is uniquely given by Figures \ref{tiling3type} and \ref{tiling3a}. The tiling is isohedral, periodic, and covers a minimal torus tiling with three tiles. 
\end{theorem}

\begin{figure}[htp]
\centering
\begin{tikzpicture}[>=latex, scale=1]

\pgfmathsetmacro{\ra}{sqrt(3)/2};

\foreach \x in {0,1}
\foreach \y in {0,1}
\fill[gray!30, shift={(0:3*\y)}, shift={(120*\x:3*\x)}]
	(180:1) -- (-60:1) -- (60:1) -- (180:1)
	;

\draw
	(0,0) -- (3,0) -- (60:3) -- (120:3) -- (0,0);
	
\foreach \x in {0,1,2}
{
\begin{scope}[shift={(60:1)}, rotate=120*\x, shift={(60:-1)}]
	
\draw[gray!70]
	(180:1) -- (-60:1) -- (60:1) -- (180:1);
	
\coordinate (A) at (0.5,\ra);
\coordinate (B) at (-1,0);
\coordinate (C) at (0.5,-\ra);

\foreach \a in {0,1,2}
{
\draw[shift={(60:1)}, rotate=120*\a, shift={(60:-1)}, red] (0,0.3) -- (A);
\draw[shift={(180:1)}, rotate=120*\a, shift={(180:-1)}, emerald] (0,0.3) -- (B);
\draw[shift={(-60:1)}, rotate=120*\a, shift={(-60:-1)}, blue] (0,0.3) -- (C);
}

\end{scope}
}

\fill (0,0.3) circle (0.05);

\draw[gray]
	(180:1) -- (-60:1) -- (60:1) -- (180:1);
	
\node at (-0.2,0) {$0$};
\node at (3.2,0) {$1$};
\node at (-1.75,2.5) {$\omega_3$};

\node[red] at (60:1.25) {$R$};
\node[emerald] at (185:1.25) {$G$};
\node[emerald] at (7:2.15) {$G'$};
\node[blue] at (-60:1.25) {$B$};
\node at (0.25,0.3) {$P$};
\node at (1.5,0.6) {$P_R$};
\node at (2,-1.1) {$P_B$};

%%%

\begin{scope}[shift={(6cm, 1cm)}]

\coordinate (R) at (0.5,\ra);
\coordinate (G) at (0.5,-\ra);

\foreach \a in {0,1,2}
{
\begin{scope}[rotate=120*\a]

\fill[yellow]
	(0:2) arc (30:90:2*\ra) -- (0,0) -- (-60:1) arc (-90:-30:2*\ra);

\draw[gray]
	(0:2) arc (30:90:2*\ra) -- (-60:1) arc (-90:-30:2*\ra);

\end{scope}
}

\node[emerald] at (2.25,0.15) {$G'$};
\node[red] at (60:0.75) {$R$};
\node[emerald] at (-1.25,0) {$G$};
\node[blue] at (0.4,-1.1) {$B$};

\begin{scope}[red]

\fill (-1,0.5) circle (0.05);

\draw (G) -- (-1,0.5) -- (R);
 
\draw[shift={(60:1)}, rotate=120, shift={(60:-1)}] (-1,0.5) -- (R);
\draw[shift={(-60:1)}, rotate=-120, shift={(-60:-1)}] (-1,0.5) -- (G);

\foreach \c in {1,-1}
\draw[shift={(0:2)}, rotate=120*\c, shift={(180:-1)}] (-1,0.5) -- (-1,0);

\end{scope}

\draw (G) -- (1,0.8) -- (R);

\fill (1,0.8) circle (0.05);
 
\draw[shift={(60:1)}, rotate=120, shift={(60:-1)}] (1,0.8) -- (R);
\draw[shift={(-60:1)}, rotate=-120, shift={(-60:-1)}] (1,0.8) -- (G);

\foreach \c in {1,-1}
\draw[shift={(0:2)}, rotate=120*\c, shift={(180:-1)}] (1,0.8) -- (-1,0);

\node at (1.1,1.05) {$P$};
\node at (0.25,1.55) {$P_R$};
\node[red] at (-1,0.7) {$P$};
\node at (2,-2.2) {$P_B$};

\end{scope}

\end{tikzpicture}
\caption{Moduli space ${\mc M}^{\text{III}}$ of type III minimal tilings.}
\label{tiling3b}
\end{figure}
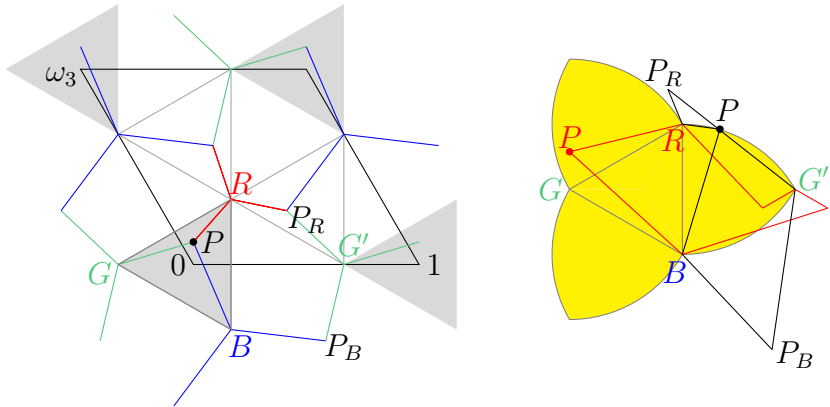

The fundamental domain of the minimal tiling and the three tiles are illustrated in the first of Figure \ref{tiling3b}. The fundamental domain can only be the specific parallelogram, and the torus can only be $T_{\omega_3}$, given by the third primitive root $\omega_3=e^{\mathrm{i}\frac{2\pi}{3}}=\frac{-1+\mathrm{i}\sqrt{3}}{2}$ of unity.

\begin{proof}
The hexagon is the tile $t_1$ in Figure \ref{tiling2type}. We indicate the sides $\bar{0}$ and $\bar{1}$ (of equal length) by red, and indicate the sides $\bar{2}$ and $\bar{3}$ by blue, and indicate the sides $\bar{4}$ and $\bar{5}$ by green. 

Consider six tiles $t_2,\dots,t_7$ around $t_1$. If the edge shared by $t_2,t_3$ is blue, then the vertex between $t_1,t_2,t_3$ is $1_12_22_3$. By the angle sum of the vertex and $[1]=\frac{2}{3}\pi$, we get $[2]=\frac{2}{3}\pi$, a contradiction. By similar reason, we know the edge between $t_2,t_3$ is not green. Therefore the edge is red, and the vertex between $t_1,t_2,t_3$ is $1_11_21_3$. By the similar argument, the vertex between $t_1,t_4,t_5$ is $3_13_43_5$, and the vertex between $t_1,t_6,t_7$ is $5_15_65_7$.

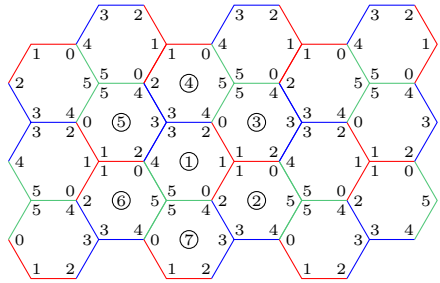
\begin{figure}[htp]
\centering
\begin{tikzpicture}[>=latex,scale=1]

\pgfmathsetmacro{\ra}{(3*sqrt(3)/10)};

\foreach \b in {-1,0,1}
{
\begin{scope}[xshift=1.8*\b cm]

\draw[red]
	(-60:0.6) -- (0:0.6) -- (60:0.6)
	(0:0.6) -- (0:1.2)
	(120:1.2) -- ++(60:0.6) -- ++(0:0.6)
	(1.2,4*\ra) -- ++(-60:0.6) -- ++(240:0.6)
	(-0.6,-2*\ra) -- ++(-60:0.6) -- ++(0:0.6);

\draw[blue]
	(180:0.6) -- (120:0.6) -- (60:0.6)
	(120:0.6) -- (120:1.2)
	(-60:0.6) -- (-60:1.2) -- ++(0:0.6)
	(-60:1.2) -- ++(240:0.6)
	(0.3,3*\ra) -- ++(60:0.6) -- ++(0:0.6)
	(1.2,2*\ra) -- ++(-60:0.6) -- ++(240:0.6);
	
\draw[emerald]
	(180:0.6) -- (240:0.6) -- (-60:0.6)
	(240:0.6) -- (240:1.2)
	(60:0.6) -- (60:1.2) -- ++(0:0.6)
	(60:1.2) -- ++(120:0.6)
	(1.2,0) -- ++(-60:0.6) -- ++(240:0.6);

\foreach \b in {0,...,5}
{
\node at (60*\b-60:0.45) {\tiny \b};
\node[yshift=2*\ra cm] at (60*\b+60:0.45) {\tiny \b};
\node[yshift=-2*\ra cm] at (60*\b+180:0.45) {\tiny \b};
\node[shift={(0.9cm, \ra cm)}] at (60*\b+180:0.45) {\tiny \b};
\node[shift={(0.9cm, 3*\ra cm)}] at (60*\b-60:0.45) {\tiny \b};
\node[shift={(0.9cm, -1*\ra cm)}] at (60*\b+60:0.45) {\tiny \b};
}

\end{scope}
}

\node[inner sep=0.5, draw, shape=circle] at (0,0) {\tiny 1};

\foreach \a in {2,...,7}
\node[inner sep=0.5, draw, shape=circle] at (60*\a-150:2*\ra) {\tiny \a};

\end{tikzpicture}
\caption{Tiling by generic hexagon of type III.}
\label{tiling3type}
\end{figure}

Since we already determined two red sides of $t_3$ and two blue sides of $t_4$, we know the edge between $t_3,t_4$ is neither red nor blue. Therefore the edge is green. This determines $t_3,t_4$. By the similar reason, we determine $t_2,t_5,t_6,t_7$. 

We see that a tile determines six tiles around it. By repeatedly applying the fact, we get the whole tiling in Figure \ref{tiling3type}. 

Let  ${\mc T}_{\text{III}}$ be the unique tiling in Proposition \ref{combo3}. It is schematically given by Figure \ref{tiling3type}, and more realistically given by Figure \ref{tiling3a}. There are red vertices, where three red edges meet. There are also the similar blue and green vertices. We obtain regular triangles by connecting vertices of three different colors. 

The symmetry $G({\mc T}_{\text{III}})$ of the tiling is the same as the symmetry of these vertices. We have the subgroup $L\cong {\bb Z}^2$ of $G({\mc T}_{\text{III}})$ consisting of the shiftings. Moreover, let $\rho$ be the rotation by $120^{\circ}$ around any colored vertex. Then we get symmetries $L\sqcup L\rho\sqcup L\rho^2$ inside $G({\mc T}_{\text{III}})$ that is transitive on all the tiles. Since the generic type III prototile is not symmetric, we get 
\[
{\mc T}_{\text{III}}
=G({\mc T}_{\text{III}})=L\sqcup L\rho\sqcup L\rho^2
=L\rtimes\langle\rho\rangle
\cong{\bb Z}^2\rtimes {\bb Z}_3.
\]

\begin{figure}[htp]
\centering
\begin{tikzpicture}[>=latex, scale=0.5]

\pgfmathsetmacro{\ra}{sqrt(3)/2};

\foreach \y in {-1,...,2}
\foreach \z in {-2,...,2}
\fill[shift={(0:3*\z)}, shift={(60:3*\y)}, gray!30]
	(180:1) -- (-60:1) -- (60:1) -- (180:1);

\foreach \x in {0,1,2}
\foreach \y in {-1,0,1}
\foreach \z in {-2,...,2}
{
\begin{scope}[shift={(0:3*\z)}, shift={(60:3*\y)}, shift={(60:1)}, rotate=120*\x, shift={(60:-1)}]

\draw[gray!70]
	(180:1) -- (-60:1) -- (60:1) -- (180:1);
	
\coordinate (A) at (0.5,\ra);
\coordinate (B) at (-1,0);
\coordinate (C) at (0.5,-\ra);

\fill[red] (A) circle (0.1);
\fill[emerald] (B) circle (0.1);
\fill[blue] (C) circle (0.1);

\foreach \a in {0,1,2}
{
\draw[shift={(60:1)}, rotate=120*\a, shift={(60:-1)}, red] (0,0.3) -- (A);
\draw[shift={(180:1)}, rotate=120*\a, shift={(180:-1)}, emerald] (0,0.3) -- (B);
\draw[shift={(-60:1)}, rotate=120*\a, shift={(-60:-1)}, blue] (0,0.3) -- (C);
}

\end{scope}
}

\end{tikzpicture}
\caption{A real example of the tiling ${\mc T}_{\text{III}}$.}
\label{tiling3a}
\end{figure}
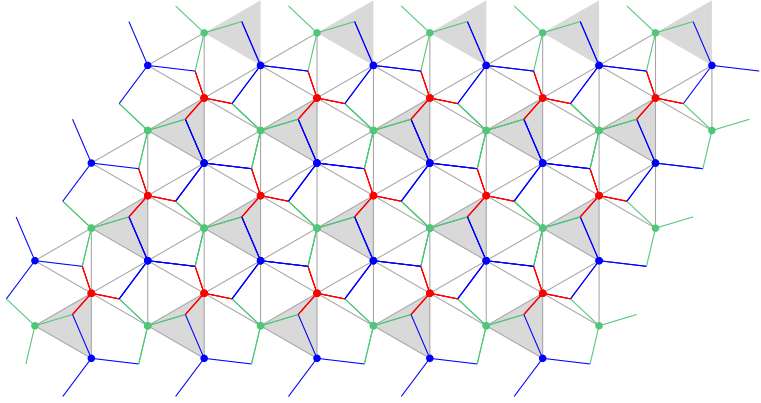

The translation subgroup $L$ is exactly all the isometries without fixed points. This implies that a type III torus tiling is ${\mc T}_{\text{III}}/G$ on ${\bb C}/G$, for a subgroup $G\sub L$ of finite index. The tiling covers the tiling ${\mc T}_{\text{III}}/L$ on the torus ${\bb C}/L$. Therefore ${\mc T}_{\text{III}}/L$, with three tiles, is the minimal type III tiling. 

Besides the colored vertices, we also have vertices where edges of three colors meet. We pick one such vertex and fix the triangle with vertices given by the other ends of the three edges. One such regular triangle and all its translations are indicated by gray in Figure \ref{tiling3a}. By putting the center of a gray triangle at the origin $0$, and set the distance between the centers of nearby gray triangles to be $1$, we identify the torus ${\bb C}/L$ with $T_{\omega_3}$, where $\omega_3=e^{\mathrm{i}\frac{2\pi}{3}}=\frac{-1+\mathrm{i}\sqrt{3}}{2}$ is the primitive cubic root of unity. 
\end{proof}

Figure \ref{tiling3b} shows the moduli space ${\mc M}^{\text{III}}$ of minimal type III tilings. In the left picture, the fundamental domain for the torus is spanned by $1$ and $\omega_3$. Around the origin is the gray triangle with vertices at $R=\frac{1}{3}e^{\mathrm{i}\frac{\pi}{3}}$, $G=\frac{1}{3}e^{\mathrm{i}\pi}$, $B=\frac{1}{3}e^{-\mathrm{i}\frac{\pi}{3}}$. We also have the point $G'$ opposite to $G$. For a point $P$, we get red edge $RP$ and blue edge $BP$. Then we rotate $RP$ around $R$ by $120^{\circ}$ to get $RP_R$, and rotate $BP$ around $B$ by $-120^{\circ}$ to get $BP_B$. Then we connect $P_R$ and $P_B$ to $G'$ to get two green edges of a type III hexagon. 

The moduli space ${\mc M}^{\text{III}}$ is the set of locations of $P$, such that the hexagon constructed above is simple. In the right of Figure \ref{tiling3b}, for the red $P$, we get the simple (and concave) red hexagon. The black $P$ is the extreme case that the (black) hexagon becomes not simple, when $P$ touches $G'P_R$. Note that $\triangle PRP_R$ is an isosceles triangle with top angle $\angle PRP_R=\frac{2}{3}\pi$. This implies $\angle G'PR=\frac{5}{6}\pi$. All the points $P$ satisfying $\angle G'PR=\frac{5}{6}\pi$ form the circular arc from $G'$ to $R$ that is centered at $B$. This is part of the boundary of the moduli space ${\mc M}^{\text{III}}$. Then we get other parts of the boundary of moduli space by symmetry.

\section{Further Discussion}

\subsection{More Examples of Minimal Tilings}

The minimal tilings for generic hexagonal prototiles of the three types consist of two, four, three tiles, respectively. There are hexagonal prototiles, such that the minimal tilings consist of a single tile. For example, the regular hexagon tiles the plane in unique way, and has a minimal tiling with a single tile on the torus $T_{\omega_3}$, $\omega_3=\frac{1+\sqrt{3}{\rm i}}{2}$. The following is a more flexible example.

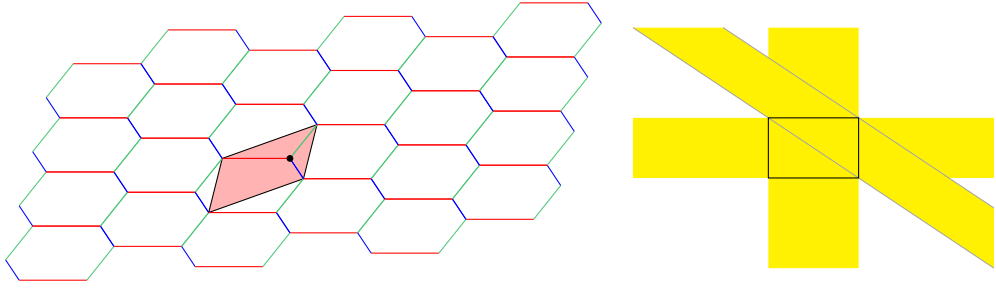
\begin{figure}[htp]
	\centering
	\begin{tikzpicture}[>=latex,scale=1]

\begin{scope}[shift={(-7.5cm,-0.5cm)}, scale=0.9]

\filldraw[fill=red!30]
	(-0.4,0.4) -- (-0.6,-0.4) -- (0.8,0.1) -- (1,0.9) -- cycle;
			
\foreach \a in {-1,0,1}
\foreach \b in {-1,...,2}
\foreach \c in {-1,1}
\foreach \d in {0,1}
{
\begin{scope}[shift={(2.6*\a cm,0.2*\a cm)}, shift={(0.2*\b cm,0.8*\b cm)}, shift={(1.4*\d cm,0.5*\d cm)}, scale=\c]	
								
\draw[red]
	(-0.4,0.4) -- (0.6,0.4);
				
\draw[emerald]
	(0.4,-0.4) -- (0.8,0.1);
				
\draw[blue]
	(0.6,0.4) -- (0.8,0.1);
								
\end{scope}
}

\fill 
	(0.6,0.4) circle (0.05);
	
\end{scope}

%% moduli

\fill[yellow]	
	(-2.4,-0.4) rectangle (2.4,0.4)
	(-0.6,-1.6) rectangle (0.6,1.6)
	(-2.4,1.6) -- (-1.2,1.6) -- (2.4,-0.8) -- (2.4,-1.6);

\draw[gray!70]
	(-1.2,1.6) -- (2.4,-0.8)
	(-2.4,1.6) -- (2.4,-1.6);
	
\draw
	(-0.6,-0.4) rectangle (0.6,0.4);
		
\end{tikzpicture}
\caption{A real example of tiling by centrally symmetric hexagon.}
\label{tiling4b}
\end{figure}

\begin{theorem}\label{combo4}
In a centrally symmetric hexagon, the opposite edges are parallel and have the same length. If the three lengths are distinct, then its tiling of the plane, such that all vertices have degree $3$, is uniquely given by Figures \ref{tiling41a} and \ref{tiling4b}. The tiling is isohedral, periodic, and covers a minimal torus tiling with a single tile. 
\end{theorem}

The fundamental domain of the minimal tiling and the single tile are illustrated in the first of Figure \ref{tiling4b}. The fundamental domain can be any parallelogram. Therefore the torus can be $T_{\tau}$ for any $\tau$ in the upper half plane.

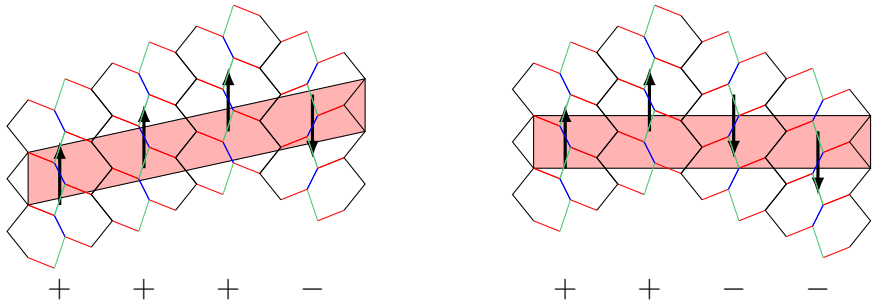
\begin{figure}[htp]
\centering
\begin{tikzpicture}[>=latex, scale=0.7]

\filldraw[fill=red!30]
	(-0.6,0.4) -- (-0.6,-0.6) -- (5.8,0.8) -- (5.8,1.8) -- cycle;
	
\filldraw[fill=red!30, xshift=9.6 cm]
	(-0.6,1.1) -- (-0.6,0.1) -- (5.8,0.1) -- (5.8,1.1) -- cycle;

\foreach \a/\b/\c in {0/0/1, 1.6/0.7/1, 3.2/1.4/1, 4.8/0.9/-1, 9.6/0.7/1, 11.2/1.4/1, 12.8/0.9/-1, 14.4/0.2/-1}
{
\begin{scope}[shift={(\a cm, \b cm)}, yscale=\c]

\draw[->, very thick]
	(0,-0.6) -- (0,0.6);
	
\foreach \d in {-1,0,1}
\foreach \z in {1,-1}
{
\begin{scope}[yshift=\d cm, scale=\z]

\draw[blue]
	(-0.1,0.2) -- (0.1,-0.2);

\draw[red]
	(0.1,0.8) -- ++(0.5,-0.2)
	(0.1,-0.2) -- ++(0.5,-0.2);

\draw[emerald]
	(0.1,0.8) -- (-0.1,0.2);
	
\draw
	(0.6,0.6) -- (1,0.1) -- (0.6,-0.4);

\end{scope}
}

\end{scope}
}	

\foreach \a in {0,1,2,6,7}
\node at (1.6*\a,-2.2) {$+$};

\foreach \a in {3,8,9}
\node at (1.6*\a,-2.2) {$-$};
			
\end{tikzpicture}
\caption{Tilings by type I hexagon satisfying $|\bar{3}|=|\bar{4}|$.}
\label{tiling4a}
\end{figure}

\begin{proof}
We indicate the three lengths by red, blue and green. In a centrally symmetric hexagon, the opposite angles are equal. Therefore we may denote the corners of a hexagon by 0, 1, 2, 0, 1, 2, as in Figure \ref{tiling4a}. Since adjacent edges have different lengths, we know the three edges at any degree $3$ vertex have different lengths. This implies $012$ is the only vertex.
	
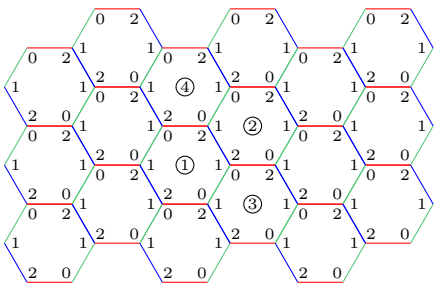
\begin{figure}[htp]
\centering
\begin{tikzpicture}[>=latex,scale=1]
			
\pgfmathsetmacro{\ra}{(3*sqrt(3)/10)};
	
\foreach \a in {-1,0,1}
\foreach \b in {1,0,-1}
\foreach \c in {1,-1}
\foreach \d in {0,1}
{
\begin{scope}[shift={(30:2*\d*\ra)}, shift={(1.8*\a cm,2*\ra*\b cm)}, scale=\c]

\draw[red]
	(120:0.6) -- (60:0.6);
\draw[emerald]
	(120:0.6) -- (180:0.6);
\draw[blue]
	(60:0.6) -- (0:0.6);

\foreach \b in {0,1,2}
\node at (60*\b-60:0.45) {\tiny \b};

\end{scope}
}

\node[inner sep=0.5, draw, shape=circle] at (0,0) {\tiny 1};
\node[inner sep=0.5, draw, shape=circle] at (-30:2*\ra) {\tiny 3};
\node[inner sep=0.5, draw, shape=circle] at (30:2*\ra) {\tiny 2};
\node[inner sep=0.5, draw, shape=circle] at (90:2*\ra) {\tiny 4};
	
\end{tikzpicture}
\caption{Tiling by centrally-symmetric hexagon.}
\label{tiling41a}
\end{figure}

We start with $t_1,t_2,t_3$ around a vertex $012$. Then the degree $3$ vertex $1_22_1\cdots=012$ determines $t_4$. More repetitions of the argument determines all the tiles.

The first of Figure \ref{tiling4b} is a real tiling of centrally-symmetric hexagon. Its symmetry group is isomorphic to ${\bb Z}^2$ and has the indicated parallelogram as the fundamental domain. Therefore the minimal tiling consists of a single tile.
\end{proof}

Given any parallelogram as the fundamental domain, the tiling is determined by the point $\bullet$. The moduli space of the minimal tilings is then all the locations of $\bullet$, such that the hexagon constructed from the point is simple. The second of Figure \ref{tiling4b} gives all such locations, with respect a rectangular fundamental domain. 

Next we show that another kind of type I hexagonal prototile can have infinitely many minimal tilings. 

\begin{theorem}\label{combo5}
A type {\rm I} hexagonal prototile satisfies 
\[
[0]+[1]+[2]=[3]+[4]+[5]=2\pi,\quad
|\bar{2}|=|\bar{5}|.
\]
Suppose it further satisfies the following conditions:
\begin{itemize}
\item Each of $\bar{0},\bar{1}$ has different lengths from the other five sides.
\item $|\bar{3}|=|\bar{4}|$.
\end{itemize}
Then its tiling of the plane, such that all vertices have degree $3$, is obtained by gluing the strips in Figure \ref{tiling1type}, and the directions of the strips can be independently chosen. Depending on the choice of directions, the tiling may or may not be periodic, and the periodic ones cover infinitely many minimal torus tilings. 
\end{theorem}

The theorem is Theorem \ref{combo1}, with the condition $|\bar{3}|\ne |\bar{4}|$ changed to $|\bar{3}|=|\bar{4}|$. This  condition is used only in the last step of the earlier proof. Before the condition is used, we already know the tiling is the union of strips. The condition $|\bar{3}|\ne |\bar{4}|$ implies all strips have the same direction. If $|\bar{3}|=|\bar{4}|$, then we can arbitrarily choose the directions of the strips. If we denote the two directions by $+$ and $-$, then a tiling corresponds to a doubly infinite sequence $S$ of the two signs.

If the sequence $S$ is not periodic, then the symmetry group is isomorphic to ${\bb Z}$, and the tiling does not cover a  torus tiling. 

If the sequence $S$ is periodic, then the symmetry group is isomorphic to ${\bb Z}^2$, and the tiling covers torus tilings. Let $p$ and $q$ be the numbers of $+$ and $-$ in one minimal period of $S$. Then the corresponding minimal tiling has $2(p+n)$ tiles. Figure \ref{tiling4a} gives two examples $+++-$ and $++--$ of the minimal period, and the corresponding fundamental domains.

Let $h$ be the minimal length of the shift that takes the strip to itself. Let $w$ be the width of the strip. Let $s$ be the offset between the two boundaries of the strip. By suitable rotation, we may assume all the strips are vertical. Then the fundamental domain is the parallelogram spanned by $(0,h)$ and $((p+n)w,(p-n)s)$. Therefore the minimal tiling is on the torus $T_{\tau}$ with  $\tau=\frac{h}{(p+n)w+\mathrm{i}(p-n)s}$.

\subsection{Tilings of a Special Torus by $12$ Tiles}

Due to the rigid nature of the fundamental domains in generic type II and type III tilings, some tori may not carry type II and type III tilings. For example, $T_{\tau}$ has type III tiling if and only if $\tau\not\in {\bb Q}+{\bb Q}\sqrt{3}\mathrm{i}$.

The torus $T_{2\sqrt{3}\mathrm{i}}$ can be tiled by generic prototiles of all three types. In Table \ref{tbe}, we list all monohedral tilings of $T_{2\sqrt{3}\mathrm{i}}$ by $12$ generic hexagons of the three types. Note that type II tilings have $l=0$ because non-zero $l$ will not give the rectangular torus $T_{\tau}$. Similarly, we cannot have $(m,n)=(2,2)$ for type III tilings. 

\begin{figure}[htp]
\centering
\begin{tikzpicture}[>=latex,scale=1]

\pgfmathsetmacro{\ra}{(sqrt(3)/3)};

%% 1

\node at (0,-0.3) {\small I(1,6;0)};
	
\begin{scope}

\foreach \b in {0,...,5}
{
\begin{scope}[xscale=0.5, yscale=0.5*\ra, yshift=1 cm+2*\b cm]

\draw
	(-0.42,1) -- (-0.5,0.8) -- (0.1,0.3) -- (-0.42,-1);

\draw[red]
	(0.6,0.1) -- (0.1,0.3)
	(-1,1) -- (-0.5,0.8);
	
\draw[blue]
	(1,-1) -- (0.6,0.1);

\draw[emerald]
	(1,1) -- (0.6,0.1);

\draw[gray!70]
	(-1,-1) rectangle (1,1);

\end{scope}
}

\end{scope}

%% 2

\node at (1.3,-0.3) {\small I(2,3;0)};

\foreach \a in {0,1}
\foreach \b in {0,1,2}
{
\begin{scope}[xshift=1cm, xscale=0.25, yscale=\ra, xshift=2*\a cm, yshift=1cm+2*\b cm]

\draw
	(-0.42,1) -- (-0.5,0.8) -- (0.1,0.3) -- (-0.42,-1);

\draw[red]
	(0.6,0.1) -- (0.1,0.3)
	(-1,1) -- (-0.5,0.8);
	
\draw[blue]
	(1,-1) -- (0.6,0.1);

\draw[emerald]
	(1,1) -- (0.6,0.1);

\draw[gray!70]
	(-1,-1) rectangle (1,1);

\end{scope}
}

%% 3

\node at (2.55,-0.3) {\small I(3,2;0)};
	
\foreach \a in {0,1,2}
\foreach \b in {0,1}
{
\begin{scope}[xshift=2.2cm, xscale=0.167, yscale=1.5*\ra, xshift=2*\a cm, yshift=1cm+2*\b cm]

\draw
	(-0.42,1) -- (-0.5,0.8) -- (0.1,0.3) -- (-0.42,-1);

\draw[red]
	(0.6,0.1) -- (0.1,0.3)
	(-1,1) -- (-0.5,0.8);
	
\draw[blue]
	(1,-1) -- (0.6,0.1);

\draw[emerald]
	(1,1) -- (0.6,0.1);

\draw[gray!70]
	(-1,-1) rectangle (1,1);

\end{scope}
}

%% 4

\node at (3.9,-0.3) {\small I(6,1;0)};

\foreach \a in {0,...,6}
\foreach \b in {0,1}
{
\begin{scope}[xshift=3.4cm, xscale=0.08333, yscale=3*\ra, xshift=2*\a cm, yshift=1cm]

\draw
	(-0.42,1) -- (-0.5,0.8) -- (0.1,0.3) -- (-0.42,-1);

\draw[red]
	(0.6,0.1) -- (0.1,0.3)
	(-1,1) -- (-0.5,0.8);
	
\draw[blue]
	(1,-1) -- (0.6,0.1);

\draw[emerald]
	(1,1) -- (0.6,0.1);

\draw[gray!70]
	(-1,-1) rectangle (1,1);

\end{scope}
}

%% 5

\node at (5.3,-0.3) {\small II(1,3;0)};

\foreach \b in {0,1,2}
{
\begin{scope}[xshift=5.3cm, xscale=0.1667, yscale=0.5*\ra, yshift=2 cm+4*\b cm]

\foreach \a in {-1,1}
{
\begin{scope}[scale=\a]

\draw
	(1,0.8) -- (-1,-0.8)
	(2.6,0.4) -- (3,0)
	(-0.4,1.6) -- (0,2)
	(3,-2) -- (2,-1.2);

\draw[red]
	(1,0.8) -- (-0.4,1.6)
	(2,-1.2) -- (3,-0.6286)
	(3,0.6286) -- (2.6,0.4);

\draw[blue]
	(1,0.8) -- (2.6,0.4)
	(0.4,-1.6) -- (2,-1.2);
	
\end{scope}
}

\draw[gray!70]
	(-3,-2) rectangle (3,2);
	
\end{scope}
}

%% 6

\node at (6.7,-0.3) {\small II(3,1;0)};

\foreach \a in {0,1,2}
{
\begin{scope}[xshift=6.3cm, xscale=0.05556, yscale=1.5*\ra, xshift=6*\a cm, yshift=2 cm]

\foreach \a in {-1,1}
{
\begin{scope}[scale=\a]

\draw
	(1,0.8) -- (-1,-0.8)
	(2.6,0.4) -- (3,0)
	(-0.4,1.6) -- (0,2)
	(3,-2) -- (2,-1.2);

\draw[red]
	(1,0.8) -- (-0.4,1.6)
	(2,-1.2) -- (3,-0.6286)
	(3,0.6286) -- (2.6,0.4);

\draw[blue]
	(1,0.8) -- (2.6,0.4)
	(0.4,-1.6) -- (2,-1.2);
	
\end{scope}
}

\draw[gray!70]
	(-3,-2) rectangle (3,2);
	
\end{scope}
}

%% 7

\node at (9,-0.3) {\small III(1,4;0)};

\begin{scope}[xshift=7.6cm, scale=0.3333]

\pgfmathsetmacro{\rra}{(sqrt(3)/2)};

\begin{scope}
\clip
	(0,0) rectangle (9, 12*\rra);

\foreach \x in {0,1,2}
\foreach \y in {0,...,4}
\foreach \z in {-2,...,3}
{
\begin{scope}[shift={(0:3*\z)}, shift={(60:3*\y)}, shift={(60:1)}, rotate=120*\x, shift={(60:-1)}]

\coordinate (A) at (0.5,\rra);
\coordinate (B) at (-1,0);
\coordinate (C) at (0.5,-\rra);

\foreach \a in {0,1}
{
\draw[shift={(60:1)}, rotate=120*\a, shift={(60:-1)}, red] (0,0.3) -- (A);
\draw[shift={(180:1)}, rotate=120*\a, shift={(180:-1)}, emerald] (0,0.3) -- (B);
\draw[shift={(-60:1)}, rotate=120*\a, shift={(-60:-1)}, blue] (0,0.3) -- (C);
}

\end{scope}
}	

\end{scope}
	
\foreach \b in {0,...,3}
\draw[shift={(120:3*\b)}, gray!40, dashed]
	(6,0) -- ++(120:3) -- ++(3,0) -- ++(-60:3) -- cycle;

\foreach \b in {0,...,3}
\draw[gray!70, yshift=3*\rra*\b cm]
	(6,0) rectangle (9,3*\rra);
	
\end{scope}

%% 8

\begin{scope}[xshift=11.6cm]

\foreach \a in {0,1,2}
\foreach \b in {0,1,2}
{
\begin{scope}[xshift=1cm, xscale=0.25, yscale=\ra, xshift=2*\a cm, yshift=1cm+2*\b cm]

\end{scope}
}

\foreach \a in {0,1,2}
\foreach \b in {0,1,2}
{
\begin{scope}[xscale=0.25, yscale=\ra, xshift=2*\a cm-0.6667*\b cm, yshift=1cm+2*\b cm, x={(1 cm,0cm)}, y={(-0.3333 cm,1 cm)}]

\draw
	(-0.42,1) -- (-0.5,0.8) -- (0.1,0.3) -- (-0.42,-1);

\draw[red]
	(0.6,0.1) -- (0.1,0.3)
	(-1,1) -- (-0.5,0.8);
	
\draw[blue]
	(1,-1) -- (0.6,0.1);

\draw[emerald]
	(1,1) -- (0.6,0.1);

\draw[gray!40, dashed]
	(-1,-1) -- (-1,1) -- (1,1) -- (1,-1) -- cycle;

\end{scope}
}

\draw[gray!70, xshift=-0.166667 cm]
	(0,0) rectangle (1,6*\ra);

\node at (0.5,-0.3) {\small I(2,3;1)};

\end{scope}
			
\end{tikzpicture}
\caption{Tilings of $T_{2\sqrt{3}i}$ with $12$ generic tiles.}
\label{tsex}
\end{figure}
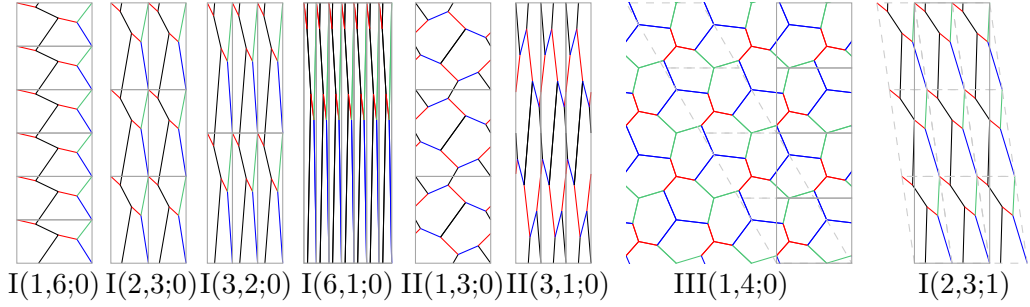

For type I, the 12-tile tiling is a six fold covering of a minimal tiling on $T_{\tau}$ consisting of two tiles. Then we have
\[
{\bb Z}+{\bb Z}2\sqrt{3}\mathrm{i}
=\lambda({\bb Z}m+{\bb Z}(l+n\tau)),\quad mn=6.
\]
For example, for $(m,n)=(2,3)$ and $l=0$, we have $\tau=\frac{2}{3}2\sqrt{3}\mathrm{i}=\frac{4}{\sqrt{3}}\mathrm{i}$. The fundamental domain for $T_{\frac{4}{\sqrt{3}}\mathrm{i}}$ is one gray rectangle in the second of Figure \ref{tsex}, with two type I tiles as indicated. The fundamental domain of the six fold cover consists of the six gray rectangles in the picture, and the covering tiling of $T_{2\sqrt{3}\mathrm{i}}$ consists of 12 tiles.
 
\begin{table}[htp]
\centering
\caption{Tilings of  $T_{2\sqrt{3}\mathrm{i}}$ with $12$ generic tiles}
\label{tbe}
\begin{tabular}{|c|c|l||c|c|l|}
\hline
Type & $(m,n)$&$l$&Type & $(m,n)$ &$l$ \\\hline
\hline
\multirow{4}{*}{I} & (1,6) & 0 
&\multirow{2}{*}{II} & (1,3) & 0 \\
& (2,3) & 0,1 & & (3,1) & 0 \\
\cline{4-6}
& (3,2) & 0,1,2 & 
\multirow{2}{*}{III} & \multirow{2}{*}{(1,4)}&\multirow{2}{*}{0} \\
& (6,1)& 0,1,2,3,4,5 & & & \\
\hline
\end{tabular}
\end{table}

We remark that the tiles may not be positioned as in the picture. We may apply actions by $SL(2,{\bb Z})$ to get other tilings as in Figure \ref{tiling1d}. Moreover, for $(m,n)=(2,3)$ and $l=1$, we get the tiling on $T_{2\sqrt{3}i}$ in the last of Figure \ref{tsex}, where the fundamental domain is given by the solid gray rectangle. 

Figure \ref{tsex} gives all tilings of $T_{2\sqrt{3}i}$ with $12$ generic tiles and satisfying $l=0$, and one tiling with $l=1$. The first four are type I tilings. Each has the moduli space with four free parameters, together with further connected components from the actions by $SL(2,{\bb Z})$. The fifth and sixth are type II tilings. Each has the moduli space with four free parameters, and no further actions by $SL(2,{\bb Z})$. The seventh is type III tiling. The moduli space has two parameters.

The special prototiles in Theorem \ref{combo4} and \ref{combo5} also give more tilings on $T_{2\sqrt{3}i}$ with 12 tiles. We choose not to pursue any more details here. A follow-up work of ours will try to identify all other special cases like in Theorem \ref{combo4} and \ref{combo5}, i.e. to get a full classification.

\subsection{Conformal Drawings of Flat Torus Tilings in $\mathbb{R}^3$}\label{sec4.1c}

Although flat tori cannot be isometrically embedded in $\mathbb{R}^3$, Pinkall \cite{pi} proved that they can be conformally embedded in $\mathbb{R}^3$ as an algebraic surface, using the Hopf fibration $\mathbb{S}^3\to\mathbb{S}^2$ and the stereographic projection of $\mathbb{S}^3$ to $\mathbb{R}^3$. In fact, the Hopf fibers over a simple closed curve $\gamma$ in $\mathbb{S}^2$ form a flat torus $T_{\frac{A+\mathrm{i} L}{4\pi}}$, where $L$ is the length of $\gamma$ and $A$ is the area enclosed by $\gamma$.

\begin{figure}[htp]
\centering
\begin{tikzpicture}[>=latex,scale=1]
	
\pgftext{
	\includegraphics[scale=0.3]{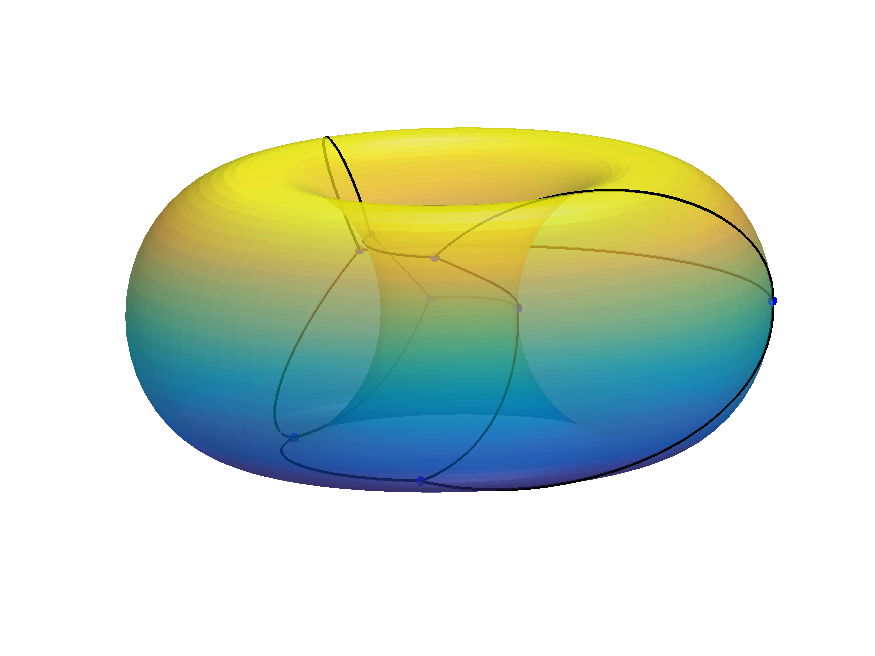}};

\node at (0,-1.6) {${\mc T}_{\rm I}(2,1;0), f=4$};
	
\begin{scope}[xshift=5cm]

\pgftext{
	\includegraphics[scale=0.3]{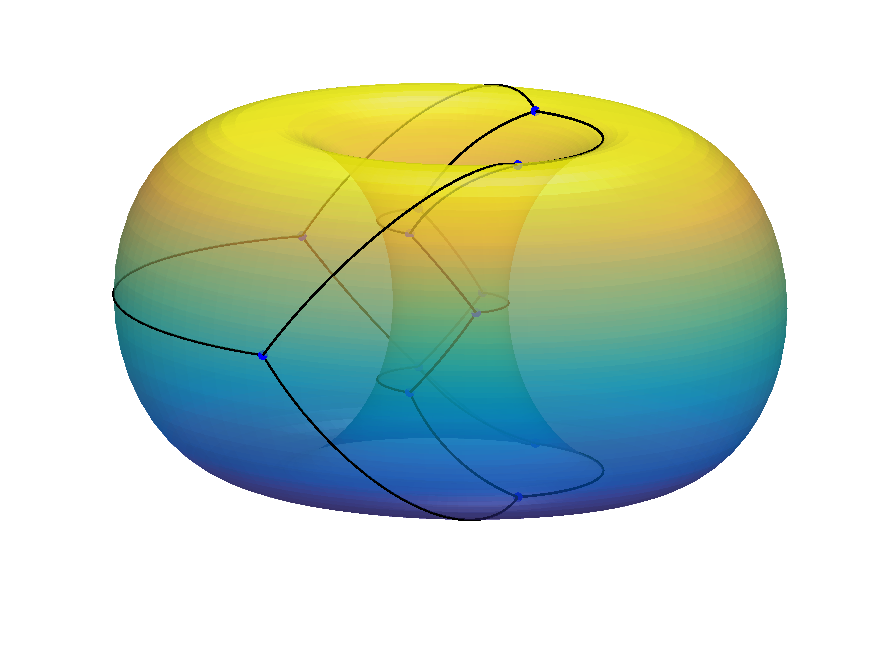}};
\node at (0,-1.6) {${\mc T}_{\rm I}(3,1;0), f=6$};

\end{scope}

\begin{scope}[xshift=10cm]

\pgftext{
	\includegraphics[scale=0.3]{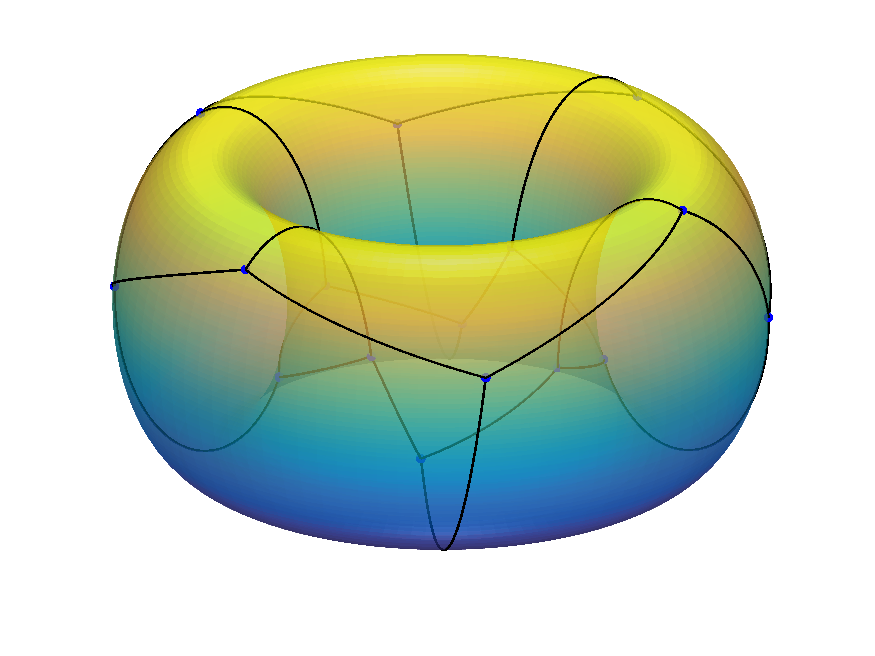}};
\node at (0,-1.6) {${\mc T}_{\rm II}(2,1;0), f=8$};

\end{scope}

\begin{scope}[yshift=-4cm]

\pgftext{
	\includegraphics[scale=0.3]{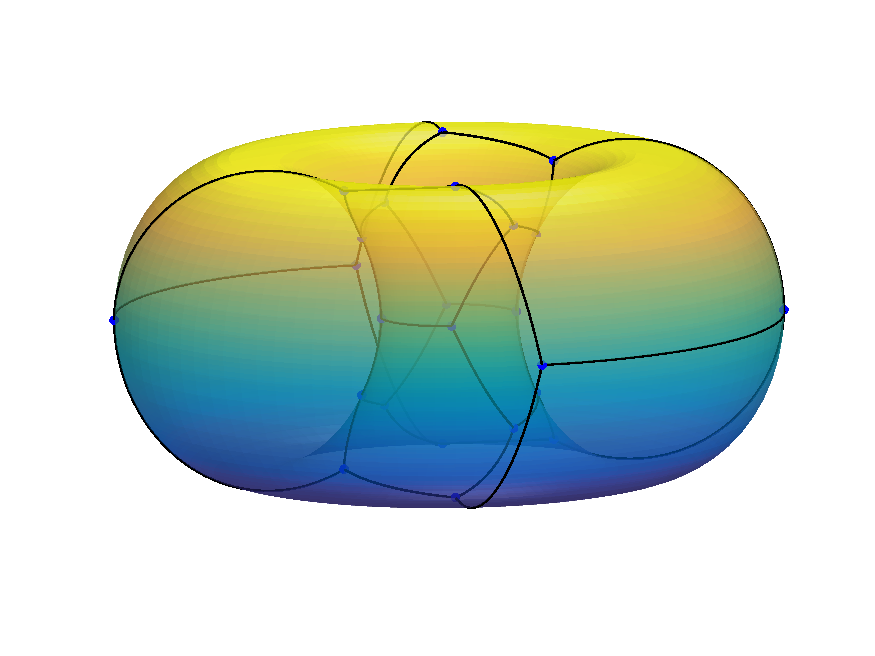}};
\node at (0,-1.6) {${\mc T}_{\rm III}(4,1;2), f=12$};
	
\begin{scope}[xshift=5cm]

\pgftext{
	\includegraphics[scale=0.3]{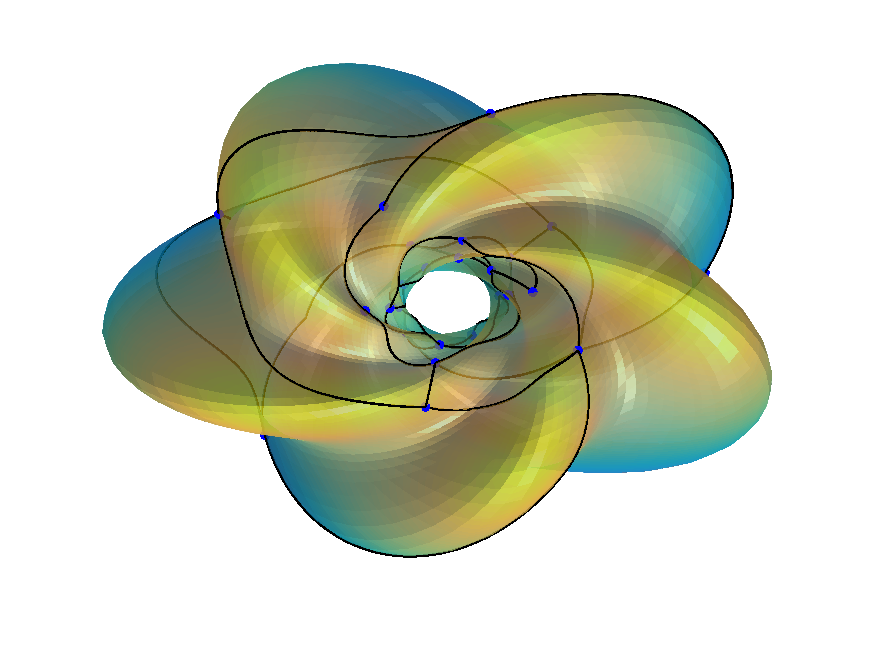}};
\node at (0,-1.6) {${\mc T}_{\rm II}(5,1;0), f=10$};

\end{scope}

\begin{scope}[xshift=10cm]

\pgftext{
	\includegraphics[scale=0.3]{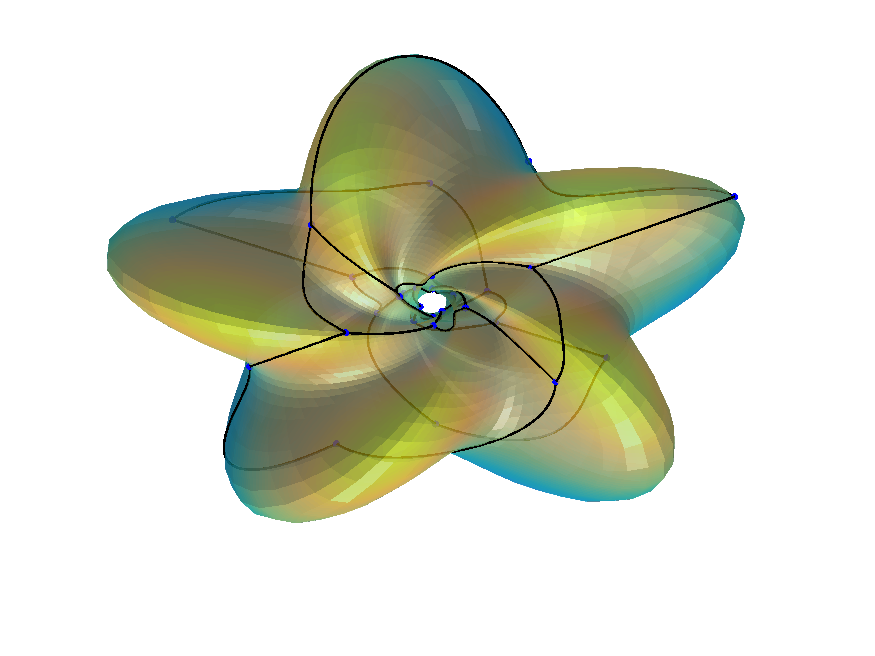}};
\node at (0,-1.6) {${\mc T}_{\rm III}(2,2;0), f=12$};
	
\end{scope}

\end{scope}
	
\end{tikzpicture}
\caption{Conformal tilings of the  torus.}
\label{61}
\end{figure}

For example, the following map conformally embeds the rectangular torus $T_{\mathrm{i} a}$ into $\mathbb{R}^3$ as a round torus
\[
x+\mathrm{i}y 
\longmapsto \frac{(a\cos2\pi x, a\sin2\pi x,a\cos \frac{2\pi y}{a})}{\sqrt{a^2+1}+\cos\frac{2\pi y}{a}}. 
\]
The torus in the first four of Figure \ref{61} are $T_{\mathrm{i} a}$ with respective $a=\frac{7}{6},1,\frac{2+\sqrt{3}}{3},\frac{2}{\sqrt{3}}$, and we draw the tilings using the embedding formula above.

For general flat tori $T_{\tau}$, where $\tau$ may not be purely imaginary, the embedding formula is much more complicated. We use explicit parametrizations in Banchoff \cite{ba}, and take the curve $\gamma$ to be
\[
\gamma(\theta)=(\sin(a+b\sin k\theta)e^{\mathrm{i}\theta}, \cos(a+b\sin k\theta)), \quad 
0 \leq \theta < 2\pi. 
\]
We use the embedding to draw the fifth and sixth tilings in Figure \ref{61}, which are on $T_{\tau}$ with respective $\tau=\frac{-1+10\mathrm{i}}{6},\frac{1+\sqrt{3}\mathrm{i}}{2}$. The tiles may twist around the meridian and longitude of a torus several times, and it takes some efforts to identify one tile. 

\end{document}